\documentclass[review,onefignum,onetabnum]{siamart171218}

\usepackage{amsfonts}
\usepackage{amssymb}
\usepackage{amsmath}
\usepackage{subfig}
\usepackage{graphicx}
\input{tcilatex}
\ifpdf
  \DeclareGraphicsExtensions{.eps,.pdf,.png,.jpg}
\else
  \DeclareGraphicsExtensions{.eps}
\fi

\headers{Convexification method for nonlinear SAR imaging}{Michael V. Klibanov et al.}

\title{Convexification Inversion Method for Nonlinear SAR Imaging with Experimentally Collected Data \thanks{Submitted to the editors DATE.
\funding{The work of Klibanov, Khoa, Smirnov, Nguyen, Bidney and Astratov was supported by US Army Research Laboratory and US Army Research Office grant W911NF-19-1-0044.}}}

\author{Michael V. Klibanov\thanks{Corresponding author. Department of Mathematics and Statistics, University of North Carolina at Charlotte, Charlotte, NC, 28223 USA  (\email{mklibanv@uncc.edu}).}
\and Vo Anh Khoa\thanks{Department of Mathematics and Statistics, University of North Carolina at Charlotte, Charlotte, NC, 28223 USA (\email{anhkhoa.vo@uncc.edu}, \email{loc.nguyen@uncc.edu}).}	
\and Alexey V. Smirnov\thanks{Department of Applied Mathematics, University of Waterloo, Waterloo, Ontario, Canada N2L 3G1 (\email{a2smirno@uwaterloo.ca}).}
\and Loc H. Nguyen\footnotemark[3]
\and Grant W. Bidney\thanks{Department of Physics and Optical Science, University of North Carolina at Charlotte, Charlotte, NC 28223, USA (\email{gbidney@uncc.edu}, \email{astratov@uncc.edu}).}
\and Lam H. Nguyen\thanks{U.S. Army Research Laboratory, Adelphi, MD 20783-1197, USA (\email{lam.h.nguyen2.civ@mail.mil}, \email{anders.j.sullivan.civ@mail.mil}).}
\and Anders J. Sullivan\footnotemark[6]
\and Vasily N. Astratov\footnotemark[5]
}

\begin{document}

\maketitle

\begin{abstract}
  This paper is concerned with the study of a version of the globally convergent convexification method with direct application to synthetic aperture radar (SAR) imaging. Results of numerial testing are presented for experimentally collected data for a fake landmine.
  The SAR imaging technique is a common tool used to create
  maps of parts of the surface of the Earth or other planets. Recently, it has been applied in the context of non-invasive inspections of buildings in military and civilian services. Nowadays, any SAR imaging software is based on the Born approximation, which is a linearization of the original wave-like partial differential equation. One of the essential assumptions this linearization procedure needs is that only those dielectric constants are imaged whose values are close to the constant background. In this work, we propose a radically new idea: to work without any linearization while still using the same data as the conventional SAR imaging technique uses. We construct a 2D image of the dielectric constant function using a number of 1D images of this function obtained via solving a 1D coefficient inverse problem (CIP) for a hyperbolic equation. Different from our previous studies on the convexification method with concentration on the global convergence of the gradient projection method, this time we prove the global convergence of the gradient descent method, which is easier to implement numerically.
\end{abstract}



\section{Introduction}

SAR imaging is a popular technique broadly used
by aircraft, drones and satellites to map some parts of the surface of
planets including Earth; cf. e.g. \cite{Amin2011,Carn1999,Gilman2017,Rotheram1985}. The SAR imaging technique
exists for several years and is developed so well by now that those flying
devices commonly use a commercial SAR software. The Physics behind it can be
understood in the following manner. It is known that the size of an antenna
of any flying device is small. However, if combining those antennas running
along a large distance covered by that device, then the size of that
synthetic antenna becomes quite large. As a result, one obtains a good
resolution of SAR images.


According to \cite{Gilman2017} and references cited therein, the
foundational assumption of any SAR imaging software is the Born
approximation. From the mathematical standpoint, this means that the
dielectric constant of the medium is assumed to be a small perturbation of
the unity. In other words, the corresponding Coefficient Inverse Problem
(CIP) for a wave-like hyperbolic partial differential equation (PDE) is
linearized near that unity. This assumption essentially enables one to
obtain an explicit formula for the solution of that CIP, which is
undoubtedly attractive from numerical standpoint. Up-to-date, any commercial
software for SAR imaging relies on that formula as well as some of its
variants in different circumstances of the above-mentioned applications. In
our first work of this direction \cite{Klibanov2021}, we have demonstrated
numerically that the Born approximation leads to significant errors in
values of the dielectric constants of targets. 


It is worth mentioning that the motivation of \cite{Klibanov2021} about nonlinear SAR imaging stems from
the non-invasive inspections of buildings with a particular concentration on
the identification of improvised explosive devices (IEDs); cf. \cite%
{Ahmad2007,Amin2011,Nguyen2008}. Needless to say, having knowledge of the
dielectric constant of the medium in that perspective is an important piece
of information which hopefully might help to decrease the false alarm rate
for military purposes. With this application being in mind, there have been
numerous models with distinctive numerical approaches tackled by the first
author and his research team for almost a decade; cf. e.g. \cite%
{Beilina2012,Khoa2020a,Khoa2020,Klibanov2019,Klibanov2015,Kolesov2017,Nguyen2018,Thanh2015}%
. The interests therein dwell not only in the approximations of the
dielectric constants, but in locations, sizes and even shapes of such
targets. This information is anticipated to be helpful in the development of
future classification algorithms that would better distinguish between
explosives and clutter in the battlefield. While the conventional SAR
imaging might only be useful in detecting the object's location and shape,
our wish here is to fill the gap regarding the identification problem of the
material of the constituent medium. Except of our work \cite{Klibanov2021}, we are unaware of any
publication, where the dielectric constants would be indeed computed for SAR
data, albeit the theory of SAR imaging tells one that those should be imaged 
\cite{Gilman2017}.

The work \cite{Khoa2020} is our first attempt to develop our convexification
method (see this section below about this method) for a CIP for the
Helmholtz equation with the moving point source. Next, we have shown
that the technique of \cite{Khoa2020} works well on experimentally collected
data \cite{Khoa2020a,Khoa2020b}. The case of the moving point source is somewhat similar to
the case of SAR imaging. We however remark that the SAR technique itself
follows time-dependent models rather than what we have done so far in \cite{Khoa2020a,Khoa2020b,Khoa2020} with the frequency domain case with a fixed frequency.
Yet, another difference between those recent works and the SAR scenario is
that in the SAR case the incident waves are sent from a specific altitude
angle and are radiated not by a point source but originated in a 2-D area, which acts as an antenna. Therefore, a \emph{radically new} idea has to be proposed to
overcome these difficulties.

We also note that various CIPs with the moving source and a fixed
frequency were considered by the group of R. G. Novikov since 1988 \cite{Novikov1989}, see \cite{Agaltsov2018,Alekseenko2008,Novikov2005,Novikov2015}. However, statements of CIPs in those publications are different
from ones in \cite{Khoa2020a,Khoa2020b,Khoa2020} and reconstruction procedures are also different
from the convexification.

The CIP of our interest with the SAR data is actually an underdetermined one 
\cite{Klibanov2021}. Indeed, the sought function that characterizes the
spatial distribution of the dielectric constant depends on three variables,
while the SAR data depemd on two variables. It is the underdertermination of
the SAR data which causes the \emph{main challenge} of the entire effort. To
address this, we form a 2D image via solutions of a number of 1D CIPs for a
wave-like hyperbolic equation. Following the approach commenced in \cite%
{Smirnov2020}, we have established in \cite{Klibanov2021} a transformation
of the underlying hyperbolic equation into a nonlocal nonlinear PDE. Then
the solution is approximated using the minimization of a weighted globally
strictly convex Tikhonov-like functional.

As to the minimization problem, we emphasize that we solve that 1D CIP by a
version of the convexification method, which is presented in \cite%
{Smirnov2020}. The convexification is a concept that the research team of
the first author has been pursuing for a number of years. This concept is
based on the Bukhgeim--Klibanov method \cite{Bukhgeim1981}, which was
originally introduced in 1981 only for proofs of global uniqueness theorems
for CIPs and without any idea at that time of its significant potential in
applications to numerical methods for CIPs. The key ingredient of the
convexification method is the construction of a Tikhonov-like cost
functional weighted by a suitable Carleman weight function. The presence of the
Carleman weight function is needed to \textquotedblleft
convexify\textquotedblright\ the cost functional involving the corresponding
PDE operator. Driven by the so-called Carleman estimate, this method is well
known to be a \emph{globally convergent} numerical method for CIPs. We call
a numerical method for a CIP globally convergent if there is a theorem which
claims that this method delivers at least one point in a sufficiently small
neighborhood of the true solution without any advanced knowledge of this
neighborhood. In the past works, the convexification method and its variants were designed
due to the need of solving CIPs that involves nonlinear differential
operators. It is obvious that for nonlinear cases, the conventional least
squares method can be trapped in local minima and ravines; cf. e.g. \cite[Figure 1]{Scales1992} for a typical graphical observation when local minima occur.
In the convexification, we are capable to prove the strict convexity of the
functional on a given convex bounded set of a Hilbert space. The global
convergence is done when one can prove the convergence of the scheme towards
the true solution starting from any element in that bounded set, whose
diameter is arbitrary.

As any CIP is highly nonlinear, we usually call our technique
\textquotedblleft nonlinear SAR imaging\textquotedblright . We point out
that our method is a heuristic one since we cannot derive it from Physics.
We can justify our idea only numerically. At the same time, we note that the
conventional SAR imaging also uses a number of heuristic approximations (in
which the Born approximation is the first one) to finalize their explicit
formula. The only rigorous part of our method is the global convergence
theory of the convexification method for our CIP.

Although the essence of this theory is taken from our recent work \cite%
{Smirnov2020}, we still have a significantly new analytical element here.
More precisely, in all previous works on the convexification (e.g., \cite%
{Khoa2020a,Khoa2020,Klibanov2018,Klibanov2019a,Klibanov2015,Smirnov2020a,Smirnov2020}%
) the global convergence was established only for the gradient projection
method, which is much harder to implement numerically than the regular
gradient descent method. On the other hand, here we establish, \emph{for the
	first time}, that actually the gradient descent method being applied to our
weighted Tikhonov-like cost functional converges globally. This result
explains our consistently repeated observation of all previous studies of
convexification, in which only the gradient descent method was numerically
implemented and numerical results were good ones; see, e.g. \cite%
{Bakushinskii2017,Khoa2020a,Khoa2020b,Khoa2020,Klibanov2019,Klibanov2019a,Klibanov2021,Smirnov2020a,Smirnov2020}%
. Hence, the extended theory in this work is also a new mathematical aspect
of our first publication \cite{Klibanov2021} which is a more engineering one. Compared with \cite{Klibanov2021}, in this paper we present  new numerical results for nonlinear SAR imaging. Numerical
results are presented for experimentally
collected SAR data. Experimental data were collected at the University of
North Carolina at Charlotte.

The paper is organized as follows. The contemporary setting of SAR imaging
is shortly introduced in section \ref{sec:2}. Then section \ref{sec:3} is
devoted to the statements of the forward and inverse problems. We also
provide an asymptotic analysis of the incident wave of the forward problem
in Lemma \ref{lem1}. This is also a novelty of the previous work \cite%
{Klibanov2021} from the theoretical standpoint. Following \cite{Smirnov2020}%
, we revisit the formulation of the 1D CIP in section \ref{sec:41}. Here, we
also present the convexification method for this CIP. Our main results are
then in section \ref{sec:43} for the global convergence of the gradient
descent method; cf. Theorems \ref{thm:11} and \ref{thm:bigthm}. As soon as
all mathematical apparatus is completed, we work in section \ref%
{sec:experiments} with our experimental results. In that
place, we prove experimentally that our proposed method delivers accurate
values of the dielectric constant of an explosive-like target. In this regard, the experimental data are
those backscattering ones for the experimental target buried
in the sandbox. Finally, some conclusions follow in section \ref%
{sec:conclusions}.

\section{Settings of SAR imaging}

\label{sec:2}

We consider in this paper the simplest version of the SAR imaging, called
stripmap SAR geometry (cf. e.g. \cite{Showmana}). In the context of the
inspection of buildings, the SAR device consists of an antenna and a
receiver. They move concurrently along a straight line and cover a certain
finite interval of this line. Due to the large antenna/targets distances, we
assume below, when discussing the inverse problem only, that the radiating
antenna and the receiver coincide and actually, this is a point. This is a
typical assumption of SAR imaging; cf. \cite{Gilman2017}. While moving, the
antenna radiates pulses of a time resolved component of the electric wave
field, and the receiver collects the backscattering time resolved signal.
Even though the propagation of the electromagnetic wave field is governed by
Maxwell's equations, the SAR imaging works with only a single wave-like PDE
governing the propagation of a single component of the electric wave field 
\cite{Gilman2017}.

Let $\Omega =(-\overline{R},\overline{R})^{3}\subset \mathbb{R}^{3}$ for $%
\overline{R}>0$ be a cube in which our targets of interest dwell. Given $%
y_{0},z_{0}\in \mathbb{R}$, we denote $\mathbf{x}_{0}=(x_{0},y_{0},z_{0})\in 
\mathbb{R}^{3}$ by the single point, where the transmitter and the receiver
are assumed to coincide. Henceforth, by letting $x_{0}$ run along an
interval $(-L,L)$ for $L>0$, we define the line of sources as follows: 
\begin{equation*}
	L_{\text{src}}:=\left\{ \mathbf{x}_{0}=\left( x_{0},y_{0},z_{0}\right) \in 
	\mathbb{R}^{3}:x_{0}\in \left( -L,L\right) \right\} ,
\end{equation*}%
where numbers $y_{0}$ and $z_{0}$ are fixed. Consider the straight line passing
through the point $\mathbf{x}_{0},$ orthogonal to the line of sources $L_{%
	\text{src}}$ and forming the angle $\theta \in \left( 0,\pi /2\right) $ with
the plane $\left\{ z=z_{0}\right\} .$ The half $CRL(\mathbf{x}_{0})$ of this
straight line with the end point $\mathbf{x}_{0},$ which looks towards the
domain of our interest $\Omega $, is called \textquotedblleft the central
line" of the antenna. The angle $\theta $ remains the same for all positions of
the source, is called the \textquotedblleft elevation angle" and it
characterizes the direction of the incident wave field. Thus, central lines
of all antennas are parallel to each other. Consider the plane $P$ which is
passing through all these central lines as well as through the line of
sources. $P$ is called \textquotedblleft slant range\textquotedblright\ or
\textquotedblleft imaging plane\textquotedblright . We remark that in the
context of SAR imaging, the domain of interest $\Omega $ should be
practically located far from the line of source $L_{\text{src}}$.

\begin{figure}[htbp]
	\centering
	\includegraphics[scale = 0.35]{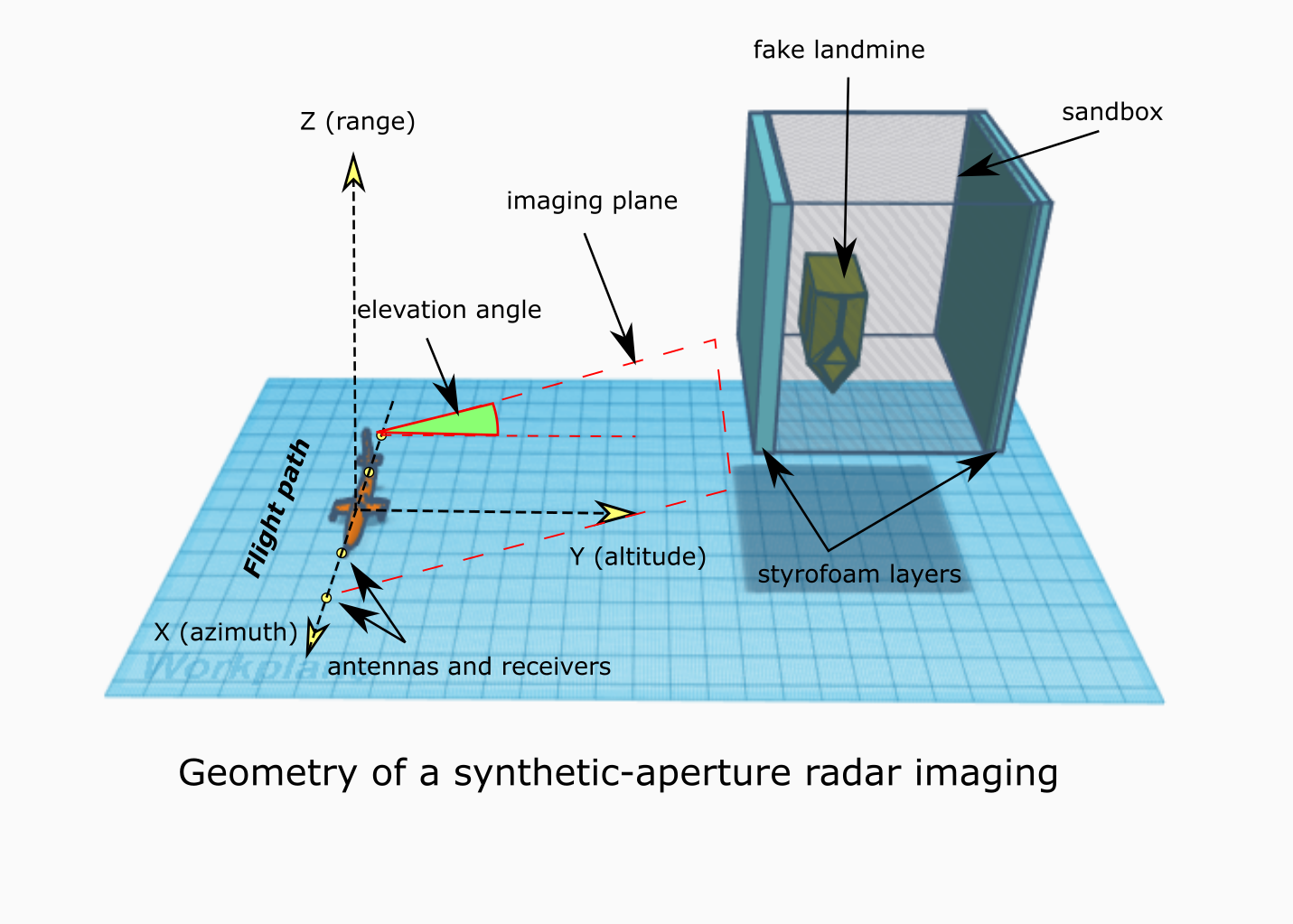}
	\caption{A schematic diagram for a stripmap SAR imaging. This depiction
		follows our experimental configuration in section \protect\ref%
		{sec:experiments}. Buried in a wooden-like box filled with the dry sand, the
		experimental target is actually a bottle of water mimicking a fake landmine. Here the transmitter and receiver run along a straight line
		in the $x-$direction. Their positions coincide and are denoted as sources $\mathbf{x}_{0}$ (yellow dots). Our
		idea is that the 2D image on the slant range (red dashed imaging plane) is
		formed by combining parallel 1D images on that plane. Those one-dimensional
		images are formed via solutions to the CIP stated in section \protect\ref%
		{sec:5.1}.}
	\label{fig:1}
\end{figure}

These settings of the geometry of the SAR imaging are depicted in
Figure \ref{fig:1}. We have prepared a glass bottle filled with clear water for our experiment; see Figure \ref{fig:5a} for some details about this
item. As to the application
being in mind, this bottle is well-suited for mimicking the usual Glassmine
43 (\cite{ChiefofOrdnance}), a non-metallic anti-personnel land mine having a glass body. It was used
in the World War II era when the Germans wanted to obstruct the detection
service of the Allies; cf. \cite{ChiefofOrdnance}. 

In order to work with the inverse problem of SAR imaging, we pose the
forward problem in the next section. This step helps us to understand the qualitative behavior of the wave field, which is incident by a close to realistic antenna, see Lemma \ref{lem1} and our comment below its proof.

\section{Statements of forward and inverse problems}

\label{sec:3}

\subsection{Forward problem} \label{sec:3.1}

For $\mathbf{x}=\left( x,y,z\right) \in \mathbb{R}^{3}$, consider a smooth
function $\varepsilon _{r}(\mathbf{x})$ that represents the spatially
distributed dielectric constant of the medium. We assume that 
\begin{equation}\label{eps}
	\begin{cases}
		\varepsilon _{r}\left( \mathbf{x}\right) \geq 1 & \text{for }\mathbf{x}\in
		\Omega , \\ 
		\varepsilon _{r}\left( \mathbf{x}\right) =1 & \text{for }\mathbf{x}\in 
		\mathbb{R}^{3}\backslash \Omega .%
	\end{cases}%
\end{equation}
In \eqref{eps}, $\varepsilon _{r}\left( \mathbf{x}\right) =1$ means that we
have vacuum outside of the domain of interest.

Now, let $T>0$ and $\mathbf{x}_{0}\in L_{\text{src}}$. Cf. \cite{Gilman2017}%
, we work below with the following initial value (or Cauchy) problem: 
\begin{equation}\label{forward1}
	\begin{cases}
		\varepsilon _{r}\left( \mathbf{x}\right) u_{tt}=\Delta u+W\left( t,\mathbf{x}%
		,\mathbf{x}_{0}\right) & \text{for }\mathbf{x}\in \mathbb{R}^{3},t\in \left(
		0,T\right) , \\ 
		u\left( \mathbf{x},\mathbf{x}_{0},0\right) =u_{t}\left( \mathbf{x},\mathbf{x}%
		_{0},0\right) =0 & \text{for }\mathbf{x}\in \mathbb{R}^{3}.%
	\end{cases}%
\end{equation}%
Here, $u=u(\mathbf{x},\mathbf{x}_{0},t)$ is the amplitude of the time
resolved component of the electric field, and the source function $W$ is
defined as:
\begin{equation}
	W\left( t,\theta ,\mathbf{x}-\mathbf{x_{0}}\right) =A\left( t\right)
	e^{-\text{i}\omega _{0}t}\widetilde{W}\left( \theta ,\mathbf{x}-\mathbf{x}%
	_{0}\right) ,\quad A(t)=\chi _{\tau }(t)e^{-\text{i}\alpha \left( t-\tau /2\right)
		^{2}},  \label{50.4}
\end{equation}%
where $\omega _{0}$ is the carrier frequency, $\alpha $ is the chirp rate,
and $\tau $ is its duration. In SAR imaging, the expression \eqref{50.4} is
called \textquotedblleft linear modulated pulse\textquotedblright\ or
\textquotedblleft chirp\textquotedblright\ with 
\begin{equation*}
	\chi _{\tau }\left( t\right) =%
	\begin{cases}
		1 & \text{for }t\in \left( 0,\tau \right) , \\ 
		0 & \text{for }t<0.%
	\end{cases}%
\end{equation*}

Different from what we have assumed above that the antenna is a point
source, when working with the data generation, we take into account a more
realistic scenario: to use a regular circular antenna. This is presented
in our establishment of the function $\widetilde{W}$ in \eqref{50.4} as
follows: The circular antenna is a disk of the diameter $D>0$ in the plane $%
P_{\text{ant}}$. Denoted by $S(\mathbf{x}_{0},\theta ,D)$, this disk is
centered at $\mathbf{x}_{0}$. Now, $\mathbf{x}_{0}$ lies in the line $L_{%
	\text{scr}}$ that intersects with the plane $P_{\text{ant}}$. On the other
hand, $P_{\text{ant}}$ is orthogonal to the central line $CRL(\mathbf{x}%
_{0}) $. Next, we change coordinates via a rotation $\mathbf{x}\rightarrow 
\mathbf{x}^{\prime }=(x^{\prime },y^{\prime },z^{\prime })$ such that the
central line $CRL(\mathbf{x}_{0})$ of the antenna is parallel to the $%
z^{\prime }$-axis. Then $\mathbf{x}_{0}=(x_{0},y_{0}^{\prime },z_{0}^{\prime
})\in L_{\text{scr}}$.

Given a number $\eta\in(0,D/2)$, we introduce a smooth function $%
m(\eta,\theta,\mathbf{x},\mathbf{x}_0,D)$: 
\begin{align}  \label{mm}
	m\left(\eta,\theta,\mathbf{x},\mathbf{x}_{0},D\right)=%
	\begin{cases}
		1 & \text{for }\left|\mathbf{x}-\mathbf{x}_{0}\right|<D/2-\eta, \\ 
		0 & \text{for }\left|\mathbf{x}-\mathbf{x}_{0}\right|\ge D/2, \\ 
		\in\left[0,1\right] & \text{for }D/2-\eta\le\left|\mathbf{x}-\mathbf{x}%
		_{0}\right|<D/2.%
	\end{cases}%
\end{align}
In this regard, $m(\eta,\theta,\mathbf{x},\mathbf{x}_0,D)=0$ outside of the
ball centered at $\mathbf{x}_0$ with the radius $D/2$. Henceforth, we define
the function $\widetilde{W}$ in \eqref{50.4} as 
\begin{align}
	\widetilde{W}(\theta,\mathbf{x}-\mathbf{x}_0) = \delta(z^{\prime }-
	z_0^{\prime })m(\eta,\theta,\mathbf{x}-\mathbf{x}_0,D),
\end{align}
where $\delta$ is the so-called delta function. By this way, due to the
presence of $m$ in \eqref{mm}, the source function $W$ vanishes outside of
our disk antenna $S(\mathbf{x}_0,\theta,D)$.

\begin{remark}
	We assume that the linear size of the antenna is small compared with the
	distance between that antenna and the domain of interest $\Omega $. This is
	equivalent with the assumption that the target is located far from the antenna.
	In physical optics, the far-field region is conceptually determined by the
	Fresnel criterion. Cf. \cite[Chapter 2]{Gilman2017}, the full angular
	beamwidth of the underlying antenna is $\Psi =2\lambda _{0}/D$, where $%
	\lambda _{0}=(2\pi c_{0})/\omega _{0}$ is the wavelength of the signal of
	the antenna and $c_{0}$ is the speed of light in the vacuum. Such
	information is only used in the SAR data generation as our best attempt to
	model the realistic scenario.
\end{remark}

The SAR data is actually the function $F(\mathbf{x}_{0},t)=u(\mathbf{x}_{0},%
\mathbf{x}_{0},t)$ for $\mathbf{x}_{0}\in L_{\text{scr}}$, $t>0$. We
generate these data via solving the forward problem \eqref{forward1}. When
doing so, we apply the Fourier transform with respect to the time variable $%
t $, which reads as 
\begin{equation}\label{fourier}
	v(\mathbf{x},\mathbf{x}_{0},k)=\int_{0}^{\infty }u(\mathbf{x},\mathbf{x}%
	_{0},t)e^{ikt}dt,
\end{equation}%
assuming that the integral converges. Then, it is straightforward to obtain
the Helmholtz equation for the resulting function $v$. Subsequently, in view
of the fact that $\varepsilon _{r}(\mathbf{x})=1$ in $\mathbb{R}%
^{3}\backslash \Omega $, we are led to an analog of the so-called
Lippmann--Schwinger equation (cf. e.g. \cite{Colton2013}): 
\begin{equation}\label{forward2}
	v\left( \mathbf{x},\mathbf{x}_{0},k\right) =v_{0}\left( \mathbf{x},\mathbf{x}%
	_{0},k\right) +k^{2}\dint\limits_{\Omega }\frac{\text{exp}\left( \text{i}%
		k\left\vert \mathbf{x}-\mathbf{x}^{\prime }\right\vert \right) }{4\pi
		\left\vert \mathbf{x}-\mathbf{x}^{\prime }\right\vert }\left( \varepsilon
	_{r}\left( \mathbf{x}^{\prime }\right) -1\right) v\left( \mathbf{x}^{\prime
	},\mathbf{x}_{0},k\right) d\mathbf{x}^{\prime }.
\end{equation}
Here, the incident wave $v_{0}\left( \mathbf{x},\mathbf{x}_{0},k\right) $
has the following form: 
\begin{align*}
	& v_{0}\left( \mathbf{x},\mathbf{x}_{0},k\right) :=\dint\limits_{0}^{\tau }A\left( t\right) e^{-\text{i}\omega _{0}t}e^{%
		\text{i}kt}dt \\
	& \times \dint\limits_{S\left( \mathbf{x}_{0},\theta ,D\right) }\frac{\text{%
			exp}\left( \text{i}k\sqrt{\left( x-\xi _{1}\right) ^{2}+\left( y-\xi
			_{2}\right) ^{2}+\left( z-z_{0}\right) ^{2}}\right) }{4\pi \sqrt{\left(
			x-\xi _{1}\right) ^{2}+\left( y-\xi _{2}\right) ^{2}+\left( z-z_{0}\right)
			^{2}}}m\left( \eta ,\theta ,\boldsymbol{\xi },\mathbf{x}_{0},D\right) d\xi
	_{1}d\xi _{2}.
\end{align*}%
In \eqref{forward2}, $\text{i}=\sqrt{-1}$. The number $k>0$ is known as the
wavenumber and we will solve \eqref{forward2} for $k\in \lbrack k_{1},k_{2}]$%
, where $k_{1}$ and $k_{2}$ are chosen numerically in our simulations. Thus,
applying the inverse Fourier transform to the function $v$ with the integration over that
interval $[k_{1},k_{2}]$ to \eqref{fourier} would result in our targeted SAR
data $F(\mathbf{x}_{0},t)$.

\begin{lemma}
	\label{lem1} The following asymptotic formula is valid: 
	\begin{align}
		& \dint\limits_{S\left( \mathbf{x}_{0},\theta ,D\right) }\frac{\text{exp}%
			\left( \text{i}k\sqrt{\left( x-\xi _{1}\right) ^{2}+\left( y-\xi _{2}\right)
				^{2}+\left( z-z_{0}\right) ^{2}}\right) }{4\pi \sqrt{\left( x-\xi
				_{1}\right) ^{2}+\left( y-\xi _{2}\right) ^{2}+\left( z-z_{0}\right) ^{2}}}%
		m\left( \eta ,\theta ,\boldsymbol{\xi },\mathbf{x}_{0},D\right) d\xi
		_{1}d\xi _{2}  \label{asymp} \\
		& =%
		\begin{cases}
			\frac{\text{i}}{2k}e^{\text{i}k\left\vert z-z_{0}\right\vert }\left[ m\left(
			\eta ,\theta ,\mathbf{x},\mathbf{x}_{0},D\right) +\mathcal{O}\left(
			1/k\right) \right] & \text{if }\left( x,y\right) \in S\left( \mathbf{x}%
			_{0},\theta ,D\right) , \\ 
			\mathcal{O}\left( 1/k^{2}\right) & \text{if }\left( x,y\right) \notin
			S\left( \mathbf{x}_{0},\theta ,D\right) .%
		\end{cases}
		\notag
	\end{align}
\end{lemma}

\textbf{Proof.}
Consider polar coordinates for the left-hand side of \eqref{asymp} with the center at $(x,y)$:
\begin{align*}
	\xi_1 - x = r\cos(\varphi),\quad
	\xi_2 - y = r\sin(\varphi).
\end{align*}
To simplify the presentation, we take $m(\xi_1,\xi_2) = m\left(\eta,\theta,\boldsymbol{\xi},\mathbf{x}_{0},D\right)$ and denote the left-hand side of \eqref{asymp} by $I_S$. Then, we arrive at
\begin{align*}
	I_S =\int_{0}^{2\pi}\int_{r_{1}\left(\varphi,x,y\right)}^{r_{2}\left(\varphi,x,y\right)}\frac{\exp\left(\text{i}k\sqrt{r^{2}+\left(z-z_{0}\right)^{2}}\right)}{4\pi\sqrt{r^{2}+\left(z-z_{0}\right)^{2}}}m\left(x+r\cos\left(\varphi\right),y+r\sin\left(\varphi\right)\right)rdrd\varphi.
\end{align*}

Notice that
\begin{align*}
	& \left(r_{1}\left(\varphi,x,y\right),\varphi\right)=\begin{cases}
		\left(0,\varphi\right) & \text{if }\left(x,y\right)\in\overline{S\left(\mathbf{x}_{0},\theta,D\right)},\\
		\in \partial S\left(\mathbf{x}_{0},\theta,D\right) & \text{if }\left(x,y\right)\notin\overline{S\left(\mathbf{x}_{0},\theta,D\right)},
	\end{cases}\\
	& \left(r_{2}\left(\varphi,x,y\right),\varphi\right)\in \partial S\left(\mathbf{x}_{0},\theta,D\right).
\end{align*}
Moreover, we have
\[
\frac{rdr}{\sqrt{r^{2}+\left(z-z_{0}\right)^{2}}}=d\left(\sqrt{r^{2}+\left(z-z_{0}\right)^{2}}\right).
\]
Therefore, by integration by parts with respect to $r$ twice, we find that
\begin{align}\label{II}
	I_{s} & =\frac{1}{\text{i}k}\left[\int_{0}^{2\pi}\frac{1}{4\pi}\exp\left(\text{i}k\sqrt{r_{2}^{2}+\left(z-z_{0}\right)^{2}}\right)m\left(x+r_{2}\cos\left(\varphi\right),y+r_{2}\sin\left(\varphi\right)\right)d\varphi\right]\\
	& -\frac{1}{\text{i}k}\left[\int_{0}^{2\pi}\frac{1}{4\pi}\exp\left(\text{i}k\sqrt{r_{1}^{2}+\left(z-z_{0}\right)^{2}}\right)m\left(x+r_{1}\cos\left(\varphi\right),y+r_{1}\sin\left(\varphi\right)\right)d\varphi\right] \nonumber \\
	& +\mathcal{O}\left(1/k^{2}\right).\nonumber
\end{align}

When $(x,y)\notin S\left(\mathbf{x}_{0},\theta,D\right)$, the first two terms in the right-hand side of \eqref{II} vanish due to the definition of $m$ in \eqref{mm}. When  $(x,y)\in S\left(\mathbf{x}_{0},\theta,D\right)$, only the first term vanishes and thus, we observe that
\begin{align*}
	& -\frac{1}{\text{i}k}\left[\int_{0}^{2\pi}\frac{1}{4\pi}\exp\left(\text{i}k\sqrt{r_{1}^{2}+\left(z-z_{0}\right)^{2}}\right)m\left(x+r_{1}\cos\left(\varphi\right),y+r_{1}\sin\left(\varphi\right)\right)d\varphi\right]\\
	& =\frac{\text{i}}{2k}e^{\text{i}k\left|z-z_{0}\right|}m\left(x,y\right).
\end{align*}
Thus, \eqref{asymp} follows. $\square $

Lemma \ref{lem1} implies that for large values of $k$ the incident wave $%
v_{0}\left( \mathbf{x},\mathbf{x}_{0},k\right) $ in \eqref{forward2} is
indeed basically contained in the cylinder whose base is the disk $S\left( 
\mathbf{x}_{0},\theta ,D\right) $. Furthermore, the magnitude of this field
changes with the distance $\left\vert z-z_{0}\right\vert $ only as $\mathcal{%
	O}(1/k)$.

\subsection{Inverse problem}

\label{sec:4}

In our inverse problem, we want to calculate the spatial distribution of the
dielectric constant $\varepsilon _{r}(\mathbf{x})$ in PDE of \eqref{forward1}%
, provided that \eqref{eps} holds and our SAR data $F(\mathbf{x}_{0},t)=u(%
\mathbf{x}_{0},\mathbf{x}_{0},t)$ are given for $\mathbf{x}_{0}\in L_{\text{%
		scr}}$, $t>0$. Recall that this is an underdetermined coefficient inverse
problem because the unknown function depends on three variables, while the
data depend only on two variables. Therefore, our goal of imaging the
function $\varepsilon _{r}(\mathbf{x})$ is extremely limited.

Instead of working on that imaging, we work on finding a certain function $%
\widetilde{\varepsilon }_{r}(\mathbf{x})$ on the slant range $P$ (i.e. for $%
\mathbf{x}\in P$). It is desirable that this function $\widetilde{%
	\varepsilon }_{r}(\mathbf{x})$ would characterize well the values of the
actual dielectric constant $\varepsilon _{r}(\mathbf{x})$. Furthermore, our
hope is that some sizes of the computed slant-range image will not be far
from the real sizes of real targets. It is also worth noting that the actual
target is not required to intersect with the slant range $P$. However, it is
still possible to obtain the above parameters of that target because the
radiation patterns of antennas might intersect with this target, thus
causing back reflected waves which arrive at the detectors.

\begin{remark}
	We call $\widetilde{\varepsilon }_{r}(\mathbf{x})$ \textquotedblleft the
	slant range image". From the Physics standpoint, the \textquotedblleft
	slant-range image\textquotedblright\ does not exist in a sense that it is
	not a physical object. In Physics, the scattering can be produced only by
	the real physical object. The mathematical inversion procedure described in
	the present paper takes into account  backscattering waves
	arriving to any detector along different directions. We use the assumption
	that the radiation pattern of the antenna has the same intensity in the
	slant-range plane as in the case when the antenna is pointed at the target.
\end{remark}

\section{Imaging via the solution of a coefficient inverse problem for a 1D
	hyperbolic PDE}

\label{sec:41}

In this part, we propose a way to form the image in the slant range $P$.
Prior to that, we want to mention that the original SAR data $F(\mathbf{x}%
_0,t) = u(\mathbf{x}_0,\mathbf{x}_0,t)$ for $\mathbf{x}_0 \in L_{\text{scr}}$%
, $t>0$ need to be preprocessed due to radar signals. We will discuss our
choice of preprocessing procedure in the numerical section.

Now, let $\mathbf{x}_{0,n}\in L_{\text{scr}}$ be the source location number $%
n=\overline{1,N}$. Let $f_{0,n}(t)$ be the preprocessed data of $F(\mathbf{x}%
_{0},t)$ at $\mathbf{x}_{0,n}$ for $n=\overline{1,N}$. Consider a
\textquotedblleft pseudo\textquotedblright\ variable $x\in \mathbb{R}$ and
let the function $c(x)\in C^{3}(\mathbb{R})$. We scale to $%
x\in (0,1)$ the length of the interval of our interest, i.e. the length of
the interval in the slant range of the central line $CRL(\mathbf{x}_{0,n})$
along which we want to calculate our function $\widetilde{\varepsilon }_{r}(%
\mathbf{x})$. Our idea below is to calculate the slant-range function $\widetilde{\varepsilon }_{r}(%
\mathbf{x})$ mentioned in section \ref{sec:4} via finding $N$ functions $c(x):=c_n(x)$ along each central line $CRL(\mathbf{x}_{0,n})$ and then forming the 2D image in the slant range plane out of these $N$ one-dimensional functions $c_{n}(x)$. This idea was
initiated in \cite{Klibanov2021}, whereas the related coefficient inverse
problem for a 1D hyperbolic PDE was considered in \cite{Smirnov2020}.

\subsection{Statement of the coefficient inverse problem for a 1D hyperbolic
	PDE}

\label{sec:5.1}

Given a number $\overline{c}>1$, we assume that 
\begin{equation}\label{cc}
	C^{3}(\mathbb{R})\ni c\left( x\right) =%
	\begin{cases}
		\in \left[ 1,\overline{c}\right]  & \text{for }x\in \left( 0,1\right) , \\ 
		1 & \text{for }x\leq 0\text{ and }x\geq 1.%
	\end{cases}%
\end{equation}%
Let $T>0$. In this scenario, our forward problem is formed as the following
initial value (or Cauchy) problem: 
\begin{equation}
	\begin{cases}
		c\left( x\right) u_{tt}=u_{xx} & \text{for }x\in \mathbb{R},t\in \left(
		0,T\right) , \\ 
		u\left( x,0\right) =0,\quad u_{t}\left( x,0\right) =\delta \left( x\right) 
		& \text{for }x\in \mathbb{R}.%
	\end{cases}
	\label{forward}
\end{equation}

Now we state our targeted inverse problem.

\textbf{Coefficient Inverse Problem (CIP)}. \emph{Determine the coefficient }%
$c\left( x\right) $ \emph{for} $x\in \left( 0,1\right) $\emph{\ in the PDE
	of \eqref{forward}, provided that the following two functions }$f_{0}\left(
t\right) $ \emph{and} $f_{1}\left( t\right) $\emph{\ are given:}%
\begin{equation}  \label{3.4}
	u\left( 0,t\right) =f_{0}\left( t\right) ,\quad u_{x}\left( 0,t\right)
	=f_{1}\left( t\right) ,\quad t\in \left( 0,T\right) .
\end{equation}

\begin{remark}
	\label{rem:41} The data at the point $x=0$ in \eqref{3.4} corresponds to
	those at the source position $\mathbf{x}_{0,n}$ in our SAR setting for each $%
	n=\overline{1,N}$. The Neumann-type data $f_{1}(t)$ in \eqref{3.4} is
	obtainable. Indeed, it was shown in \cite{Smirnov2020} (and also commenced
	in \cite{Sochacki1988}) that the function $u(x,t)$ in the PDE \eqref{forward} satisfies the so-called absorbing boundary conditions: 
	\begin{align}
		& u_{x}\left( a_{1},t\right) -u_{t}\left( a_{1},t\right) =0,\quad a_{1}\leq
		0,  \label{ab1} \\
		& u_{x}\left( a_{2},t\right) +u_{t}\left( a_{2},t\right) =0,\quad a_{2}\geq 
		\sqrt{\overline{c}}.  \label{ab2}
	\end{align}%
	In \eqref{ab1}, take $a_{1}=0$. Then, using \eqref{3.4}, we obtain $%
	f_{1}(t)=f_{0}^{\prime }(t)$. Thus, for each $n=\overline{1,N}$ we find the
	Neumann data $f_{1}(t)$ for our CIP by differentiating the preprocessed data 
	$f_{0,n}(t)$.
\end{remark}

\begin{remark}
	After solving the aforementioned CIP, we scale back the \textquotedblleft
	pseudo\textquotedblright\ interval $x\in (0,1)$ and obtain $x_{n}^{\prime
	}\in (0,a)$ for a certain number $a>0$. In fact, for each $n=\overline{1,N}$ we assign 
	\begin{equation*}
		\left. \widetilde{\varepsilon }_{r}\left( \mathbf{x}\right) \right\vert _{%
			\mathbf{x}\in CRL\left( \mathbf{x}_{0,n}\right) }:=\widetilde{\varepsilon }%
		_{r,n}\left( \mathbf{x}\right) :=c\left( x_{n}^{\prime }\right) ,\quad
		x_{n}^{\prime }\in \left( 0,a\right) .
	\end{equation*}%
	We define the function $\widetilde{\varepsilon }_{r}\left( \mathbf{x}%
	\right) $ as a set of values of the function $\widetilde{\varepsilon }_{r}(%
	\mathbf{x})$ on 1D cross-sections of the slant range $P$ by central lines $%
	CRL(\mathbf{x}_{0,n})$. Thus, our strategy is to form the 2D image on the
	slant range $P$ as the collection of 1D cross-sections: 
	\begin{equation*}
		\widetilde{\varepsilon }_{r}\left( \mathbf{x}\right) =\left\{ \widetilde{%
			\varepsilon }_{r,n}\left( \mathbf{x}\right) ,n=\overline{1,N}\right\} .
	\end{equation*}%
	This finalizes our steps of the SAR imaging procedure.
\end{remark}

Now we want to discuss the analysis of our CIP. When doing so, we revisit a
transformation investigated in \cite{Smirnov2020} to delineate an essential
nonlocal nonlinear PDE for our mathematical analysis.

Cf. \cite{Romanov2019}, consider the one-to-one change of variables 
\begin{equation*}
	y=y(x)=\int_{0}^{x}\sqrt{c(s)}ds.
\end{equation*}%
Then $y(x)$ is called the travel time, which the wave needs to travel from
the source $\left\{ 0\right\} $ to a certain point $\left\{ x\right\} $. By %
\eqref{cc}, $dy/dx=\sqrt{c(y)}\geq 1$. Moreover, $y(x)=x$ for $x\leq 0$ and $%
x(y)=y$ for $y\leq 0$. We introduce the following functions $w,Q,p$: 
\begin{align}
	& w\left( y,t\right) =u\left( x\left( y\right) ,t\right) c^{1/4}\left(
	x\left( y\right) \right) ,  \label{w} \\
	& Q\left( y\right) =c^{-1/4}\left( x\left( y\right) \right) ,\quad p\left(
	y\right) =\frac{Q^{\prime \prime }\left( y\right) }{Q\left( y\right) }-2%
	\left[ \frac{Q^{\prime }\left( y\right) }{Q\left( y\right) }\right] ^{2}.
	\label{Qp}
\end{align}%
Due to \eqref{cc}, functions $Q$ and $p$ are well-defined and we have 
\begin{align}
	& w(y,t)=u(x(y),t)\quad \text{for }y<0\;\text{and }y>\sqrt{\overline{c}}, 
	\notag \\
	& p\left( y\right) \in C^{1}\left( \mathbb{R}\right) ,\quad p\left( y\right)
	=0\quad \text{for }y<0\;\text{and }y>\sqrt{\overline{c}}.  \label{pp1}
\end{align}%
Therefore, we obtain from \eqref{forward}--\eqref{3.4}, \eqref{w}, \eqref{Qp}
that 
\begin{align}
	& w_{tt}=w_{yy}+p\left( y\right) w\quad \text{for }y\in \mathbb{R},t\in
	\left( 0,\widetilde{T}\right) ,  \label{w1} \\
	& w\left( y,0\right) =0,\quad w_{t}\left( y,0\right) =\delta \left( y\right)
	\quad \text{for }y\in \mathbb{R},  \label{w2} \\
	& w\left( 0,t\right) =f_{0}\left( t\right) ,\quad w_{y}\left( 0,t\right)
	=f_{1}\left( t\right) \quad \text{for }t\in \left( 0,\widetilde{T}\right) ,
	\label{w3}
\end{align}%
where the number $\widetilde{T}\geq 2\sqrt{\overline{c}}$ depends on $T$.

\begin{remark}
	It was shown in \cite[Chapter 2]{Romanov2019} and \cite[Section 2]%
	{Smirnov2020a} that the forward problem \eqref{w1}--\eqref{w2} has a unique
	solution in $C^{3}(\mathbb{R})$. In particular, it holds that 
	\begin{align}
		& w\left( y,t\right) =\frac{H\left( t-\left\vert y\right\vert \right) }{2}+%
		\frac{1}{2}\int_{\left( y-t\right) /2}^{\left( y+t\right) /2}p\left( \xi
		\right) \left( \int_{\left\vert \xi \right\vert }^{t-\left\vert y-\xi
			\right\vert }w\left( \xi ,\tau \right) d\tau \right) d\xi ,  \label{w4} \\
		& w\left( y,t\right) =0\quad \text{for }t<\left\vert y\right\vert ,\quad
		\lim_{t\rightarrow \left\vert y\right\vert ^{+}}w\left( y,t\right) =\frac{1}{%
			2},  \label{w5}
	\end{align}%
	where $H(z)$ for $z\in \mathbb{R}$ is the Heaviside function. Here, $w\in
	C^{3}(t\geq \left\vert y\right\vert )$ is guaranteed by $p(y)\in C^{1}(%
	\mathbb{R})$, which explains the condition $c(x)\in C^{3}(\mathbb{R})$
	imposed in \eqref{cc}.
\end{remark}

\begin{remark}
	\label{rem4} It is straightforward that to find $c(x)\in C^{3}(\mathbb{R})$
	stated in our CIP, we need to determine the function $p(y)\in C^{1}(\mathbb{R%
	}),$ which is the solution of the CIP \eqref{w1}--\eqref{w3}. Since $y\leq 
	\sqrt{\overline{c}}$, we find $p(y)$ for $y\in (0,b)$ for an arbitrary $%
	b\geq \sqrt{\overline{c}}$. Then the \textquotedblleft
	threshold\textquotedblright\ $\widetilde{T}$ for the CIP \eqref{w1}--\eqref{w3} is taken by a certain $\widetilde{T}\geq \max \left\{
	2b,T\right\} \geq 2\sqrt{\overline{c}}$ in order to guarantee the uniqueness
	of the CIP \eqref{w1} and \eqref{w3}; cf. \cite{Romanov2019}.
\end{remark}

Let the parameter $\beta \in (0,1/2).$ We will choose it later. Consider an
arbitrary number $\mu \in (0,2\beta b)$. We introduce the following
rectangle $R$ and the triangle $R_{\beta ,\mu ,b}\subset R$, as follows: 
\begin{align}
	& R:=\left\{ \left( y,t\right) \in \mathbb{R}^{2}:y\in \left( 0,b\right)
	,t\in \left( 0,\widetilde{T}\right) \right\} ,  \label{R1} \\
	& R_{\beta ,\mu ,b}:=\left\{ \left( x,t\right) \in \mathbb{R}^{2}:x+\beta
	t<2\beta b-\mu \;\text{for }x,t>0\right\} .  \label{R2}
\end{align}%
Henceforth, we consider a new function $v(y,t)$, 
\begin{equation}\label{vv}
	v\left( y,t\right) =w\left( y,t+y\right) ,\quad x,t>0.
\end{equation}%
Then $v\in C^{3}(y\geq 0,t\geq 0)$ and \eqref{w4}--\eqref{w5} imply 
\begin{equation}\label{ww}
	\lim_{t\rightarrow 0^{+}}v\left( y,t\right) =\lim_{t\rightarrow
		0^{+}}w\left( y,t+y\right) =\lim_{t\rightarrow y^{+}}w\left( y,t\right) =%
	\frac{1}{2}.
\end{equation}

It follows from \eqref{w1} and \eqref{ww} that 
\begin{align}
	& v_{yy}-2v_{yt}+p\left( y\right) v=0\quad \text{for }\left( y,t\right) \in
	R,  \label{v1} \\
	& v\left( y,0\right) =\frac{1}{2}\quad \text{for }y\in \left( 0,b\right) .
	\label{v2}
\end{align}%
Take $t=0$ in \eqref{v1} and use \eqref{v2}, \eqref{ww}. We deduce that 
\begin{equation}\label{pp}
	p(y)=4v_{yt}(y,0),\quad y\in (0,b).
\end{equation}

Next, we differentiate both sides of \eqref{v1} with respect to $t$ and
consider $V(y,t)=v_{t}(y,t)$ satisfying $V\in C^{2}(\overline{R})$. Using %
\eqref{3.4}, \eqref{ab2} and \eqref{vv}--\eqref{pp}, we obtain the following
nonlocal nonlinear boundary value problem for a hyperbolic equation: 
\begin{align}
	& V_{yy}-2V_{yt}+4V_{y}\left( y,0\right) V=0\quad \text{for }\left(
	y,t\right) \in R,  \label{V1} \\
	& V\left( 0,t\right) =q_{0}\left( t\right) :=f_{0}^{\prime }\left( t\right)
	\quad \text{for }t\in \left( 0,\widetilde{T}\right) ,  \label{V2} \\
	& V_{y}\left( 0,t\right) =q_{1}\left( t\right) :=f_{0}^{\prime \prime
	}\left( t\right) +f_{1}^{\prime }\left( t\right) \quad \text{for }t\in
	\left( 0,\widetilde{T}\right) ,  \label{V3} \\
	& V_{y}\left( b,t\right) =0\quad \text{for }t\in \left( 0,\widetilde{T}%
	\right) .  \label{V4}
\end{align}

If $V(y,t)$ for $(y,t)\in R$ is found via solving the overdetermined
boundary value problem \eqref{V1}--\eqref{V4}, then we can rely on \eqref{pp}
to have 
\begin{equation}\label{pp2}
	p(y)=4V_{y}(y,0),\quad y\in (0,b).
\end{equation}%
Then recalling Remark \ref{rem4} we are able to determine the coefficient $%
c(x)$ of the original CIP; see \eqref{Qp}. Thus, we need to solve the
nonlinear overdetermined boundary value problem \eqref{V1}--\eqref{V4}.

\subsection{Convexification method for the nonlocal nonlinear problem}

\label{sec:5.2}

Consider the rectangle $R$ in \eqref{R1}. Whenever we deal below with
function spaces $H^{k}(R)$, we work only with real-valued functions.
Consider now the subspaces $H_{0}^{2}(R)\subset H^{2}(R)$ and $%
H_{0}^{4}(R)\subset H^{4}(R)$, 
\begin{align*}
	& H_{0}^{2}\left( R\right) :=\left\{ u\in H^{2}\left( R\right) :u\left(
	0,t\right) =u_{y}\left( 0,t\right) ,u_{y}\left( b,t\right) =0\right\} , \\
	& H_{0}^{4}\left( R\right) :=H^{4}\left( R\right) \cap H_{0}^{2}\left(
	R\right) ,
\end{align*}%
see Remark \ref{rem:47} about $H^{4}(R).$ Recall that the parameters $\beta
\in (0,1/2)$, $b\geq \sqrt{\overline{c}}$ are introduced in \eqref{R2}. As
introduced in \cite{Smirnov2020a,Smirnov2020}, we consider $\psi _{\lambda
}(y,t)$, 
\begin{equation}\label{CWF}
	\psi _{\lambda }(y,t)=e^{-2\lambda (y+\beta t)},\quad \lambda \geq 1.
\end{equation}%
This is the Carleman Weight Function (CWF) for the linear operator $\partial
_{y}^{2}-2\partial _{y}\partial _{t}$. This operator is the principal part
of the operator of equation \eqref{V1}, which we want to solve.

\begin{theorem}[Carleman estimate \protect\cite{Smirnov2020a}]
	\label{carleman} There exist constants $C=C(\beta ,R)>0$ and $\lambda
	_{0}=\lambda _{0}(\alpha ,R)\geq 1$ depending only on listed parameters such
	that for all functions $u\in H_{0}^{2}(R)$ and for all $\lambda \geq \lambda
	_{0}$ the following Carleman estimate is valid: 
	\begin{align*}
		& \int\limits_{R}\left( u_{yy}-2u_{yt}\right) ^{2}\psi _{\lambda }dydt\geq
		C\int\limits_{R}\left( \lambda \left( u_{y}^{2}+u_{t}^{2}\right) +\lambda
		^{3}u^{2}\right) \psi _{\lambda }dydt \\
		& +C\int\limits_{0}^{b}\left( \lambda u_{y}^{2}+\lambda ^{3}u^{2}\right)
		\left( y,0\right) e^{-2\lambda y}dy-Ce^{-2\lambda \widetilde{T}%
		}\int\limits_{0}^{b}\left( \lambda u_{y}^{2}+\lambda ^{3}u^{2}\right) \left(
		y,\widetilde{T}\right) dy.
	\end{align*}
\end{theorem}

Now we introduce a Tikhonov-like cost functional $J_{\lambda ,\gamma }$
weighted by the CWF \eqref{CWF}, where $\gamma \in (0,1)$ is the
regularization parameter. Let $M>0$ be an arbitrary number. Consider the
convex set $B(M,q_{0},q_{1})\subset H^{4}\left( R\right) $ of the diameter $%
2M$, 
\begin{align}
	& B\left( M,q_{0},q_{1}\right)  \label{ball} \\
	& =\left\{ u\in H^{4}\left( R\right) :u\left( 0,t\right) =q_{0}\left(
	t\right) ,u_{y}\left( 0,t\right) =q_{1}\left( t\right) ,u_{y}\left(
	b,t\right) =0,\left\Vert u\right\Vert _{H^{4}\left( R\right) }<M\right\} . 
	\notag
\end{align}%
In \eqref{ball}, $q_{0}$ and $q_{1}$ are given data in \eqref{V2} and %
\eqref{V3}, respectively.

\begin{remark}\label{rem:47}
	We need to use the space $H^{4}(R)$ since $H^{4}(R)\subset C^{2}(\overline{R%
	})$ is known by the embedding theorem. Theorem \ref{sconvex} is a direct
	analog of Theorem 2 of \cite{Smirnov2020}, which requires the regularity $%
	C^{2}(\overline{R})$; see the inequality (73) in \cite{Smirnov2020}.
	However, in computations we use $H^{2}(R)$ instead of $H^{4}(R)$, 
	see section \ref{sec:experiments}. Below, $B\left( M,q_{0},q_{1}\right) $ denotes only
	interior points of the set $\overline{B\left( M,q_{0},q_{1}\right) }$, which
	is closed in the space $H^{4}(R)$.
\end{remark}

Following \eqref{V1}, denote 
\begin{equation}\label{SS}
	\mathcal{S}(V)=V_{yy}-2V_{yt}+4V_{y}(y,0)V,\quad (y,t)\in R.
\end{equation}%
Henceforth, our weighted Tikhonov-like cost functional $J_{\lambda ,\gamma
}:H^{5}(R)\rightarrow \mathbb{R}_{+}$ is defined as follows: 
\begin{equation}\label{JJ}
	J_{\lambda ,\gamma }\left( V\right) =\dint\limits_{R}\left[ \mathcal{S}%
	\left( V\right) \right] ^{2}\psi _{\lambda }dydt+\gamma \left\Vert
	V\right\Vert _{H^{4}\left( R\right) }^{2},\quad \gamma \in (0,1).
\end{equation}

\textbf{Minimization Problem.} \emph{Minimize the Tikhonov-like cost
	functional }$J_{\lambda ,\gamma }\left( V\right) $\emph{\ in \eqref{JJ} on
	the set }$\overline{B\left( M,q_{0},q_{1}\right) }.$

\begin{theorem}[Global strict convexity \protect\cite{Smirnov2020}]
	\label{sconvex} For all $\lambda ,\gamma >0$ and for all elements $V\in 
	\overline{B\left( M,q_{0},q_{1}\right) }$, the cost functional (\ref{JJ})
	has the Fr\'{e}chet derivative $J_{\lambda ,\gamma }^{\prime }\left(
	V\right) \in H_{0}^{4}\left( R\right) .$ Let $\lambda _{0}=\lambda
	_{0}(\beta ,R)\geq 1$ be the number obtained in Theorem \ref{carleman}. Let $%
	R_{\alpha ,\mu ,b}$ be the triangle defined in (\ref{R2}) and let $%
	\widetilde{T}\geq 2b.$ Then there exists a sufficiently large number $%
	\lambda _{1}=\lambda _{1}\left( R,M,\beta \right) \geq $ $\lambda _{0}$
	depending only on $R,M$ and $\beta $ such that for all $\lambda \geq \lambda
	_{1}$ and for all $\gamma \in \lbrack 2e^{-\lambda \beta \widetilde{T}},1),$
	the functional $J_{\lambda ,\gamma }(V)$ is strictly convex on the set $%
	\overline{B\left( M,q_{0},q_{1}\right) },$ i.e. for all $V_{1},V_{2}\in 
	\overline{B\left( M,q_{0},q_{1}\right) }$ the following estimate holds\label%
	{convex} 
	\begin{align}
		& J_{\lambda ,\gamma }\left( V_{2}\right) -J_{\lambda ,\gamma }\left(
		V_{1}\right) -J_{\lambda ,\gamma }^{\prime }\left( V_{1}\right) \left(
		V_{2}-V_{1}\right) \geq Ke^{-2\lambda \left( 2\beta b-\mu \right)
		}\left\Vert V_{2}-V_{1}\right\Vert _{H^{1}\left( R_{\beta ,\mu ,b}\right)
		}^{2} \\
		& +Ke^{-2\lambda \left( 2\beta b-\mu \right) }\left\Vert V_{2}\left(
		y,0\right) -V_{1}\left( y,0\right) \right\Vert _{H^{1}\left( 0,2\beta b-\mu
			\right) }^{2}+\frac{\gamma }{2}\left\Vert V_{2}-V_{1}\right\Vert
		_{H^{4}\left( R\right) }^{2},  \notag
	\end{align}%
	where the constant $K=K\left( \beta ,\mu ,b,M\right) >0$ depends only on
	listed parameters.
\end{theorem}

To this end, $K>0$ denotes different constants depending on some parameters,
to be specified. Theorem \ref{sconvex} claims the global strict convexity
since smallness conditions are not imposed on the diameter $2M>0$ of the
convex set $B(M,q_0,q_1)$. Below we formulate the theorem about the
existence and uniqueness of a minimizer of the functional $%
J_{\lambda,\gamma}(V)$ on the set $\overline{B(M,q_0,q_1)}$.

\begin{theorem}[Existence and uniqueness of a minimizer \protect\cite%
	{Smirnov2020}]
	\label{Vunique} Let parameters $\lambda _{1}$ and $\gamma $ be the same as
	in Theorem \ref{sconvex}. Then for any $\lambda \geq \lambda _{1}$ there
	exists a unique minimizer $V_{\min ,\lambda ,\gamma }$ of the functional $%
	J_{\lambda ,\gamma }\left( V\right) $ on the set $\overline{B\left(
		M,q_{0},q_{1}\right) }.$ Furthermore, the following inequality holds 
	\begin{equation}
		J_{\lambda ,\gamma }^{\prime }\left( V_{\min ,\lambda ,\gamma }\right)
		\left( V-V_{\min ,\lambda ,\gamma }\right) \geq 0\quad\text{for any } V\in 
		\overline{B\left( M,q_{0},q_{1}\right) }.  \label{3.20}
	\end{equation}
\end{theorem}

Using the standard concept of the regularization theory (cf. \cite%
{Beilina2012,Tikhonov1995}), we assume the existence of ideal, noiseless
data $f_{0}^{\ast }\left( t\right) ,f_{1}^{\ast }\left( t\right) $ in (\ref%
{V3}), which correspond to the exact coefficient $p^{\ast }\left( y\right)
\in C^{1}\left( \mathbb{R}\right) $ and the exact function $w^{\ast }\left(
y,t\right) $ in (\ref{w1})--(\ref{w3}). In addition, we assume that the
function $p^{\ast }\left( y\right) $ satisfies conditions (\ref{pp1}). Let
the function $V^{\ast }\left( y,t\right) $ be obtained from the function $%
w^{\ast }\left( y,t\right) $ by the same way we have obtained above $V$ from 
$w$; cf. section \ref{sec:5.1}. Assume that 
\begin{align}  \label{Vstart}
	V^{\ast}\in B(M,q_0^{\ast},q_1^{\ast}),
\end{align}
where functions $q_0^{\ast},q_1^{\ast}$ are the noiseless data of $q_0,q_1$
in \eqref{V2}, \eqref{V3}, respectively. By \eqref{V1}--\eqref{V4} and %
\eqref{SS}, we have 
\begin{align*}
	\mathcal{S}(V^{\ast})=0,\quad (y,t)\in R,\;\text{and } p^{\ast}(y) =
	4v_{y}^{\ast}(y,0),\quad y\in (0,b).
\end{align*}

Let $\sigma >0$ be a sufficiently small number which we will choose later.
This number characterizes the level of noise in the data $q_{0},q_{1}.$
Assume that $\sigma \in \left( 0,\min \left\{ 1,M\right\} \right) $. Suppose
that there exist functions $F\in B\left( M,q_{0},q_{1}\right) $ and $F^{\ast
}\in B\left( M,q_{0}^{\ast },q_{1}^{\ast }\right) $ such that 
\begin{equation}
	\left\Vert F-F^{\ast }\right\Vert _{H^{4}(R)}<\sigma .  \label{3.23}
\end{equation}

The presence of functions $F$ and $F^{\ast }$ is to obtain the zero
Dirichlet and Neumann boundary condition at $\left\{ x=0\right\} $ in %
\eqref{V2} and \eqref{V3}. In fact, for every $V\in B\left(
M,q_{0},q_{1}\right) $ we consider 
\begin{equation*}
	W^{\ast }=V^{\ast }-F^{\ast },\quad W=V-F.
\end{equation*}%
Then by the triangle inequality, it is clear that 
\begin{equation*}
	W,W^{\ast }\in B_{0}\left( 2M\right) :=\left\{ W\in H_{0}^{4}\left( R\right)
	:\left\Vert W\right\Vert _{H^{4}\left( R\right) }<2M\right\} .
\end{equation*}%
Moreover, for any $W\in B_{0}\left( 2M\right) $ it holds that $W+F\in
B\left( 3M,q_{0},q_{1}\right) $. Therefore, we can consider a modification $%
\tilde{J}_{\lambda ,\gamma }:B_{0}(2M)\rightarrow \mathbb{R}_{+}$ of the
functional $J_{\lambda ,\gamma }$, 
\begin{equation*}
	\tilde{J}_{\lambda ,\gamma }(W)=J_{\lambda ,\gamma }(W+F)\quad \text{for any 
	}W\in B_{0}(2M).
\end{equation*}%
As inherited from properties of $J_{\lambda ,\gamma }$, the functional $%
\tilde{J}_{\lambda ,\gamma }$ is globally strict convex on the ball $%
\overline{B_{0}(2M)}\subset H_{0}^{4}(R)$ and admits a unique minimizer $%
W_{\min ,\lambda ,\gamma }$ on that ball. For brevity, we do not formulate
the theorem for $\tilde{J}_{\lambda ,\gamma }$ as it is analogous to the
statement in Theorem \ref{sconvex}. The presence of the functional $\tilde{J}%
_{\lambda ,\gamma }$ is necessary to estimate the accuracy of the minimizer $%
V_{\min ,\lambda ,\gamma }$ in the next theorem; see \cite[Theorem 5]%
{Smirnov2020} for details of the proof.

\begin{theorem}[Accuracy estimate of the minimizer \protect\cite{Smirnov2020}%
	]
	\label{thm:10} Suppose that conditions (\ref{Vstart}) and (\ref{3.23}) hold.
	Assume that $\widetilde{T}\geq 4b.$ Denote 
	\begin{equation*}
		\omega =\frac{\alpha \left( \widetilde{T}-4b\right) +\mu }{2\left( 2\beta
			b-\mu \right) },\quad \rho =\frac{1}{2}\min \left\{ \omega ,1\right\} .
	\end{equation*}%
	Let $\lambda _{1}=\lambda _{1}\left( R,M,\beta \right) $ be the number of
	Theorem \ref{sconvex}. Denote $\lambda _{2}=\lambda _{1}\left( R,3M,\beta
	\right) \geq \lambda _{1}.$ Let the number $\sigma _{0}\in \left( 0,1\right) 
	$ be sufficiently small such that $\ln \sigma _{0}^{-1/\left( 2\left( 2\beta
		b-\mu \right) \right) }\geq \lambda _{2}.$ Let the number $\sigma \in \left(
	0,\sigma _{0}\right) .$ Let the numbers $\lambda =\lambda \left( \sigma
	\right) $ and $\gamma =\gamma \left( \sigma \right) $ be defined as: 
	\begin{eqnarray}
		&&\lambda =\lambda \left( \sigma \right) =\ln \sigma ^{-1/\left( 2\left(
			2\beta b-\mu \right) \right) }>\lambda _{2},  \label{3.24} \\
		&&\gamma =\gamma \left( \sigma \right) =2e^{-\lambda \beta \widetilde{T}%
		}=2\sigma ^{\left( \beta \widetilde{T}\right) /\left( 2\left( 2\beta b-\mu
			\right) \right) }.  \label{3.25}
	\end{eqnarray}%
	Let $V_{\min ,\lambda ,\gamma }\left( y,t\right) $ be the unique minimizer
	of the cost functional $J_{\lambda ,\gamma }\left( V\right) $ on the set $%
	B\left( M,q_{0},q_{1}\right) $, cf. Theorem \ref{Vunique}. Let $p_{\min
		,\lambda ,\gamma }\left( y\right) $ be the coefficient obtained from the
	minimizer $V_{\min ,\lambda ,\gamma }\left( y,t\right) $ via (\ref{pp2}),
	i.e. $p_{\min ,\lambda ,\gamma }\left( y\right) =4\partial _{y}V_{\min
		,\lambda ,\gamma }\left( y,0\right) $. Then the following accuracy estimates
	are valid 
	\begin{eqnarray}
		&&\left\Vert V_{\min ,\lambda ,\gamma }-V^{\ast }\right\Vert _{H^{1}\left(
			R_{\beta ,\mu ,b}\right) }\leq K\sigma ^{\rho },  \label{3.26} \\
		&&\left\Vert p_{\min ,\lambda ,\gamma }-p^{\ast }\right\Vert _{L^{2}\left(
			0,2\beta b-\mu \right) }\leq K\sigma ^{\rho }.  \label{3.27}
	\end{eqnarray}
\end{theorem}

\begin{theorem}[Lipschitz continuity of $J_{\protect\lambda ,\protect\gamma %
	}^{\prime }$]
	\label{lipsta} As in Theorem \ref{sconvex}, let $J_{\lambda ,\gamma
	}^{\prime }\left( v\right) \in H_{0}^{4}\left( R\right) $ be the Fr\'{e}chet
	derivative of the functional $J_{\lambda ,\gamma }\left( V\right) $. Then
	for any $\lambda ,\gamma >0$ it satisfies the Lipschitz continuity condition
	on the set $\overline{B\left( M,q_{0},q_{1}\right) }.$ More precisely, there
	exists a number $\widetilde{C}=\widetilde{C}\left( \lambda ,\gamma ,M\right)
	>0$ depending only on listed parameters such that 
	\begin{equation*}
		\left\Vert J_{\lambda ,\gamma }^{\prime }\left( V_{1}\right) -J_{\lambda
			,\gamma }^{\prime }\left( V_{2}\right) \right\Vert _{H^{4}\left( R\right)
		}\leq \widetilde{C}\left\Vert V_{1}-V_{2}\right\Vert _{H^{4}\left( R\right) }%
		\text{ for all }V_{1},V_{2}\in \overline{B\left( M,q_{0},q_{1}\right) }.
	\end{equation*}
\end{theorem}

The proof of this theorem is omitted here because it is completely
similar with the proof of Theorem 3.1 of \cite{Bakushinskii2017}.

\begin{remark}
	In our first paper on SAR imaging \cite{Klibanov2021}, we have seen a large
	amount of numerical and experimental results using the gradient descent
	method, although the convergence theory is based on the gradient projection
	method in \cite{Smirnov2020}. The only reason is that the gradient descent
	method is rather easy-to-implement numerically. Therefore, this is the
	moment to explore the convergence analysis of the gradient descent method.
	It is also accentuated that starting from the work \cite{Bakushinskii2017},
	only theorems about the global convergence of the gradient projection method
	in the context of the convexification method were studied; cf. e.g. \cite%
	{Khoa2020a,Khoa2020b,Khoa2020,Klibanov2019,Klibanov2019a,Smirnov2020a,Smirnov2020}%
	. Hence, our analysis below is really important to fulfill the gap between
	the theory and numerics of the convexification method.
\end{remark}

\subsection{Global convergence of the gradient descent method}

\label{sec:43}

Let the number $\varkappa \in (0,1)$ and the number $\vartheta \in (0,1/3)$.
Let $V^{(0)}\in B\left( \vartheta M,q_{0},q_{1}\right) $ be an arbitrary
point in that set. We set the gradient descent method as follows: 
\begin{equation}
	V^{\left( n\right) }=V^{\left( n-1\right) }-\varkappa J_{\lambda ,\gamma
	}^{\prime }\left( V^{\left( n-1\right) }\right) ,\quad n=1,2,\ldots
	\label{gradient}
\end{equation}%
Note that $V^{(n)}$ denotes the $n$th iteration that expects to be close to
the minimizer $V_{\min ,\lambda ,\gamma }$ as $n$ tends to infinity. It is
also worth noting that since $J_{\lambda ,\gamma }^{\prime }\left(
V^{(n-1)}\right) \in H_{0}^{4}(R)$ by Theorem \ref{sconvex}, then all
functions $V^{(n)}$ of the sequence \eqref{gradient} satisfy the same
boundary conditions as ones in \eqref{ball}.

As assumed in \eqref{Vstart}, we know that the exact solution $V^{\ast }$ is
an interior point of the set $B(M,q_{0}^{\ast },q_{1}^{\ast })$. This is an
interior point in terms of the norm of the space $H^{4}(R)$. On the other
hand, by \eqref{3.26} the distance between the exact solution $V^{\ast }$
and the minimizer $V_{\min ,\lambda ,\gamma }$ is small in the norm of the
space $H^{1}\left( R_{\beta ,\mu ,b}\right) $, as long as the noise level $%
\sigma $ in the data is sufficiently small. Therefore, it is reasonable to
assume that the minimizer $V_{\min ,\lambda ,\gamma }$ is an interior point
of the set $B(M,q_{0},q_{1})$ rather than the one located on the boundary of
this set, see Theorem \ref{Vunique}. More precisely, we need to assume that 
\begin{equation}
	V_{\min ,\lambda ,\gamma }\in B\left( \vartheta M,q_{0},q_{1}\right) .
	\label{newVstar}
\end{equation}

Now we are in position to estimate the distance between the approximation $%
V^{(n)}$ and the minimizer $V_{\min ,\lambda ,\gamma }$.

\begin{theorem}[Global convergence of the gradient method]
	\label{thm:11} Fix a number $\vartheta \in (0,1/3)$. Assume that %
	\eqref{newVstar} holds with this $\vartheta .$ Suppose that the parameters $%
	\lambda $ and $\gamma $ are the same as in Theorem \ref{sconvex}. Then there
	exists a sufficiently small number $\varkappa _{0}\in (0,1)$ such that the
	sequence $\left\{ V^{(n)}\right\} _{n=0}^{\infty }\subset B\left(
	M,q_{0},q_{1}\right) $ for each $\varkappa \in (0,\varkappa _{0})$.
	Furthermore, for each $\varkappa \in (0,\varkappa _{0})$ there exists a
	number $\theta =\theta (\varkappa )\in (0,1)$ such that 
	\begin{equation}\label{conv}
		\left\Vert V_{\min ,\lambda ,\gamma }-V^{\left( n\right) }\right\Vert
		_{H^{4}\left( R\right) }\leq \theta ^{n}\left\Vert V_{\min ,\lambda ,\gamma
		}-V^{\left( 0\right) }\right\Vert _{H^{4}\left( R\right) },\quad n=1,2,\ldots
	\end{equation}
\end{theorem}

\textbf{Proof.} Consider the operator $Z:B(M,q_{0},q_{1})\rightarrow
H^{4}(R) $ defined as 
\begin{equation}\label{Z}
	Z(V)=V-\varkappa J_{\lambda ,\gamma }^{\prime }(V).
\end{equation}%
It follows from the proof of \cite[Theorem 2.1]{Bakushinskii2017} that
Theorems \ref{sconvex} and \ref{lipsta} imply that there indeed exist
numbers $\varkappa _{0}\in (0,1)$ and $\theta (\varkappa )\in (0,1)$ for $%
\varkappa \in (0,\varkappa _{0})$ such that the operator $Z$ defined in %
\eqref{Z} has the contractual mapping property, which reads as 
\begin{equation}
	\left\Vert Z\left( V_{1}\right) -Z\left( V_{2}\right) \right\Vert
	_{H^{4}\left( R\right) }\leq \theta \left\Vert V_{1}-V_{2}\right\Vert
	_{H^{4}\left( R\right) }\text{ for all }V_{1},V_{2}\in B\left(
	M,q_{0},q_{1}\right) .  \label{ZZ}
\end{equation}

Since $B\left( \vartheta M,q_{0},q_{1}\right) \subset B\left(
M,q_{0},q_{1}\right) $, then estimate \eqref{ZZ} holds for all $%
V_{1},V_{2}\in B\left( \vartheta M,q_{0},q_{1}\right) $. However, it was not
proven in \cite{Bakushinskii2017} that the operator $Z$ maps the set $%
B\left( M,q_{0},q_{1}\right) $ into itself. Note that since $J_{\lambda
	,\gamma }^{\prime }\in H_{0}^{4}(R)$, then functions $Z(V)$ and $V$ satisfy
the same boundary conditions for any $V\in B\left( M,q_{0},q_{1}\right) $.

Since $V_{\min ,\lambda ,\gamma }$ is an interior point of the set $B\left(
\vartheta M,q_{0},q_{1}\right) ,$ then $J_{\lambda ,\gamma }^{\prime }\left(
V_{\min ,\lambda ,\gamma }\right) =0$. Therefore, 
\begin{equation*}
	V_{\min ,\lambda ,\gamma }=V_{\min ,\lambda ,\gamma }-\varkappa J_{\lambda
		,\gamma }^{\prime }\left( V_{\min ,\lambda ,\gamma }\right) =Z\left( V_{\min
		,\lambda ,\gamma }\right) \in B\left( \vartheta M,q_{0},q_{1}\right) .
\end{equation*}%
It now remains to check that $V^{(n)}\in B\left( M,q_{0},q_{1}\right) $
whenever $V^{(0)}\in B\left( \vartheta M,q_{0},q_{1}\right) $ for any $%
\vartheta \in (0,1/3)$. In fact, it can be done by mathematical induction. For $n=1$, we particularly have
\begin{align}\label{VVV}
	V^{(1)} = V^{(0)} - \varkappa J_{\lambda
		,\gamma }^{\prime }(V^{(0)}) = Z(V^{(0)}).
\end{align}
Therefore, using \eqref{ZZ}, we estimate that
\begin{align}\label{WWW}
	\left\Vert V^{\left( 1\right) }-V_{\min ,\lambda ,\gamma }\right\Vert
	_{H^{4}\left( R\right) }\le \theta\left\Vert V^{\left( 0\right) }
	- V_{\min ,\lambda ,\gamma } \right\Vert _{H^{4}\left(
		R\right) }.
\end{align}
Recall that $V^{(0)}\in B\left( \vartheta M,q_{0},q_{1}\right) $. By the triangle inequality, we obtain
\begin{align*}
	\left\Vert V^{\left( 1\right) }\right\Vert _{H^{4}\left( R\right) }& \leq
	\left\Vert V_{\min ,\lambda ,\gamma }\right\Vert _{H^{4}\left( R\right)
	}+\theta\left\Vert V^{\left( 0\right) }-V_{\min ,\lambda ,\gamma
	}\right\Vert _{H^{4}\left( R\right) } \\
	& \leq \left( 1+\theta\right) \left\Vert V_{\min ,\lambda ,\gamma
	}\right\Vert _{H^{4}\left( R\right) }+\theta \left\Vert V^{\left(
		0\right) }\right\Vert_{H^{4}\left( R\right) } \\
	& \leq 2\vartheta M+\theta \vartheta M<3\vartheta M<M.\text{ }
\end{align*}%
Henceforth, $V^{(1)}\in B(M,q_0,q_1)$. In the same vein as \eqref{VVV}, we have for $n=2$ that $V^{(2)} = Z(V^{(1)})$. Therefore, using \eqref{ZZ} and \eqref{WWW}, we estimate that
\begin{align}\label{WWWW}
	\left\Vert V^{\left( 2\right) }-V_{\min ,\lambda ,\gamma }\right\Vert
	_{H^{4}\left( R\right) }\le \theta\left\Vert V^{\left( 1\right) }
	- V_{\min ,\lambda ,\gamma } \right\Vert _{H^{4}\left(
		R\right) } \le \theta^2\left\Vert V^{\left( 0\right) }
	- V_{\min ,\lambda ,\gamma } \right\Vert _{H^{4}\left(
	R\right) }.
\end{align}
Applying again the triangle inequality and using \eqref{WWWW}, we obtain
\begin{align*}
	\left\Vert V^{\left( 2\right) }\right\Vert _{H^{4}\left( R\right) }& \leq
	\left\Vert V_{\min ,\lambda ,\gamma }\right\Vert _{H^{4}\left( R\right)
	}+\theta^2\left\Vert V^{\left( 0\right) }-V_{\min ,\lambda ,\gamma
	}\right\Vert _{H^{4}\left( R\right) } \\
	& \leq \left( 1+\theta^2\right) \left\Vert V_{\min ,\lambda ,\gamma
	}\right\Vert _{H^{4}\left( R\right) }+\theta^2 \left\Vert V^{\left(
		0\right) }\right\Vert_{H^{4}\left( R\right) } \\
	& \leq 2\vartheta M+\theta^2 \vartheta M<3\vartheta M<M.\text{ }
\end{align*}%
Henceforth, $V^{(2)}\in B(M,q_0,q_1)$.

By continuing this process up to the $n$-th iteration, we obtain the same that $V^{(n)}\in B(M,q_0,q_1)$ and \eqref{conv} holds true. Hence, we complete the proof of the theorem. $\square $

In the following, we prove that the sequence \eqref{gradient} of the
gradient descent method actually converges to the exact solution, as long as
the noise level $\sigma $ in the data (cf. \eqref{3.23}) tends to zero. Then
the convergence of the sequence $\left\{ p^{(n)}\right\} _{n=0}^{\infty }$,
where $p^{(n)}(y)=4\partial _{y}V_{y}^{(n)}(y,0)$ (see \eqref{pp2}), follows
immediately.

\begin{theorem}
	\label{thm:bigthm} Suppose that conditions of Theorems \ref{thm:10} and \ref%
	{thm:11} hold true. Then the following convergence estimate is valid: 
	\begin{align*}
		\left\Vert V^{\ast}-V^{\left(n\right)}\right\Vert
		_{H^{1}\left(R_{\alpha,\mu,b}\right)} & +\left\Vert
		p^{\ast}-p^{\left(n\right)}\right\Vert _{L^{2}\left(0,2\alpha b-\mu\right)}
		\\
		& \le2K\sigma^{\rho}+4C\theta^{n}\left\Vert
		V_{\min,\lambda,\gamma}-V^{\left(0\right)}\right\Vert
		_{H^{4}\left(R\right)},\quad n=1,2,\ldots,
	\end{align*}
	where $C=C(R)>0$ depends only on the rectangle $R$.
\end{theorem}

\textbf{Proof.} Combine \eqref{3.26} and \eqref{conv}. Since $%
H^{4}(R)\subset H^{1}\left( R_{\beta ,\mu ,b}\right) $, then we use the
triangle inequality, 
\begin{align*}
	\left\Vert V^{\ast }-V^{\left( n\right) }\right\Vert _{H^{1}\left( R_{\beta
			,\mu ,b}\right) }& \leq \left\Vert V^{\ast }-V_{\min ,\lambda ,\gamma
	}\right\Vert _{H^{1}\left( R_{\beta ,\mu ,b}\right) }+\left\Vert V_{\min
		,\lambda ,\gamma }-V^{\left( n\right) }\right\Vert _{H^{1}\left( R_{\beta
			,\mu ,b}\right) } \\
	& \leq K\sigma ^{\rho }+\theta ^{n}\left\Vert V_{\min ,\lambda ,\gamma
	}-V^{\left( 0\right) }\right\Vert _{H^{4}\left( R\right) }.
\end{align*}

We estimate the distance between $p^{\ast }$ and $p^{(n)}$ similarly.
Indeed, we have 
\begin{align*}
	\left\Vert p^{\ast }-p^{\left( n\right) }\right\Vert _{L^{2}\left( 0,2\beta
		b-\mu \right) }& \leq \left\Vert p^{\ast }-p_{\min ,\lambda ,\gamma
	}\right\Vert _{L^{2}\left( 0,2\beta b-\mu \right) }+\left\Vert p_{\min
		,\lambda ,\gamma }-p^{\left( n\right) }\right\Vert _{L^{2}\left( 0,2\beta
		b-\mu \right) } \\
	& \leq K\sigma ^{\rho }+4\left\Vert \partial _{y}V_{\min ,\lambda ,\gamma
	}\left( \cdot ,0\right) -\partial _{y}V^{\left( n\right) }\left( \cdot
	,0\right) \right\Vert _{L^{2}\left( 0,2\alpha b-\mu \right) } \\
	& \leq K\sigma ^{\rho }+4C\theta ^{n}\left\Vert V_{\min ,\lambda ,\gamma
	}-V^{\left( 0\right) }\right\Vert _{H^{4}\left( \Omega \right) },
\end{align*}%
where we have used \eqref{3.27}, \eqref{conv} and the trace inequality%
\begin{equation*}
	\left\Vert \mathcal{U}_{y}\left( \cdot ,0\right) \right\Vert _{L^{2}\left(
		0,2b-\mu \right) }\leq C\left\Vert \mathcal{U}\right\Vert _{H^{4}\left(
		\Omega \right) }\text{ for any } \mathcal{U}\in H^{4}(R).
\end{equation*}%

This completes the proof of the theorem. $\square $ 

\begin{remark}
	Theorem \ref{thm:11} is our central theorem since this is the first time we
	prove the global convergence of the gradient descent method by restricting
	the minimizer in a smaller set of $B\left(M,q_0,q_1\right)$. It is worth
	mentioning that the convergence in Theorem \ref{thm:11} as well as its
	consequence in Theorem \ref{thm:bigthm} are still global since the starting
	point $V^{(0)}$ is arbitrary in $B\left(\vartheta M,q_0,q_1\right)$ with $%
	\vartheta$ being fixed in $(0,1/3)$, and smallness conditions are not
	imposed on the number $M$.
\end{remark}

\section{Numerical verification using experimental data}

\label{sec:experiments}

In our first paper on SAR imaging \cite{Klibanov2021}, we verified the proposed method using the simulated data. This means the SAR data were computed by solving the Lippmann--Schwinger
equation \eqref{forward2} obtained after numerous assumptions about the
forward problem. This equation can be solved using a trigonometric
Galerkin-based numerical approach postulated in \cite{Lechleiter2013}. We also tested in \cite{Klibanov2021} our method using experimental data for a
two-sided scanning of a barrack building in which experimental objects with
high dielectric constants were located.

In this work, we want to corroborate the method using experimental data collected at the University of North Carolina at Charlotte. We have placed horizontally the watered bottle inside
of a wooden framed box fully filled with the dry sand. The front and back
sides of that sandbox were covered by styrofoam layers. The thickness of
each layer is about 7 centimeters (cm), while the sandbox itself has the
dimensions of approximately $50\times 50\times 70$ $\text{cm}^{3}$ (H$\times 
$W$\times $L). That styrofoam has almost the same dielectric constant as the
air and thus, it does not affect our backscattering signals. Following the
purpose of the search of anti-personnel landmines, the bottle is buried close to
the front surface with 10 cm depth. Getting a landmine with the 10 cm depth
is usually considered in many models of de-mining operations; cf. e.g. \cite%
{Daniels2006}. 



\begin{figure}[tbp]
	\begin{centering}
		\includegraphics[scale=0.07]{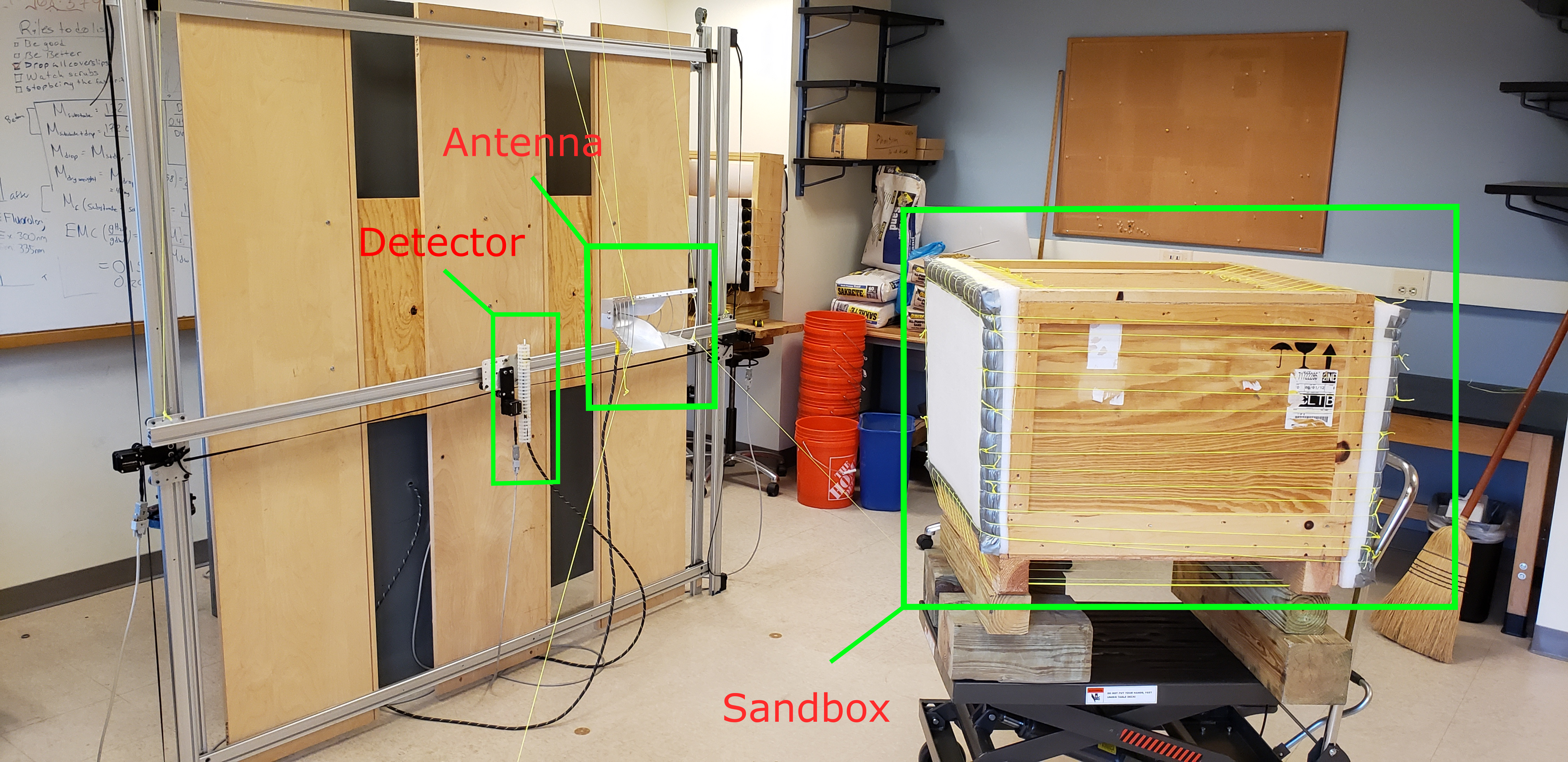}
		\par\end{centering}
	\caption{A photograph of our experimental configuration with a fake landmine
		buried in a sandbox. That sandbox was moved along a straight line which has
		formed the $\pi/4$ angle with the wall which can be seen in the
		photograph. The total length of this movement was 140 cm with the 5 cm step
		size. We have used the
		concurrent movement along a straight line of both the antenna and detector, see Figure \protect\ref{fig:1}. This models well the
		SAR data collection of Figure \protect\ref{fig:1}. }
	\label{fig:setup}
\end{figure}

\begin{figure}[tbp]
	\begin{centering}
		\subfloat[\label{fig:jug2}]{\includegraphics[scale=0.45]{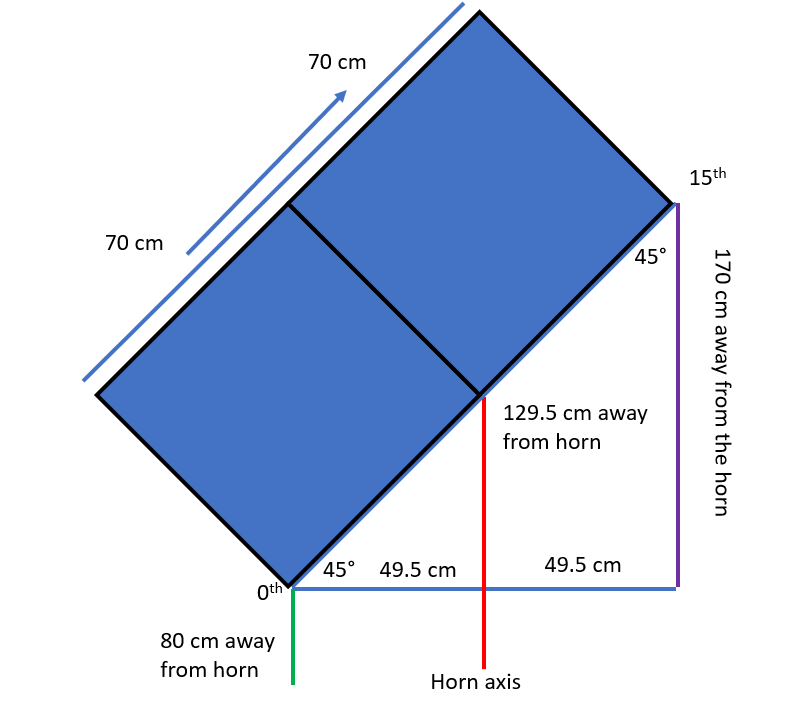}
		}
		\par\end{centering}
	\begin{centering}
		\subfloat[\label{fig:jug3}]{\includegraphics[scale=0.35]{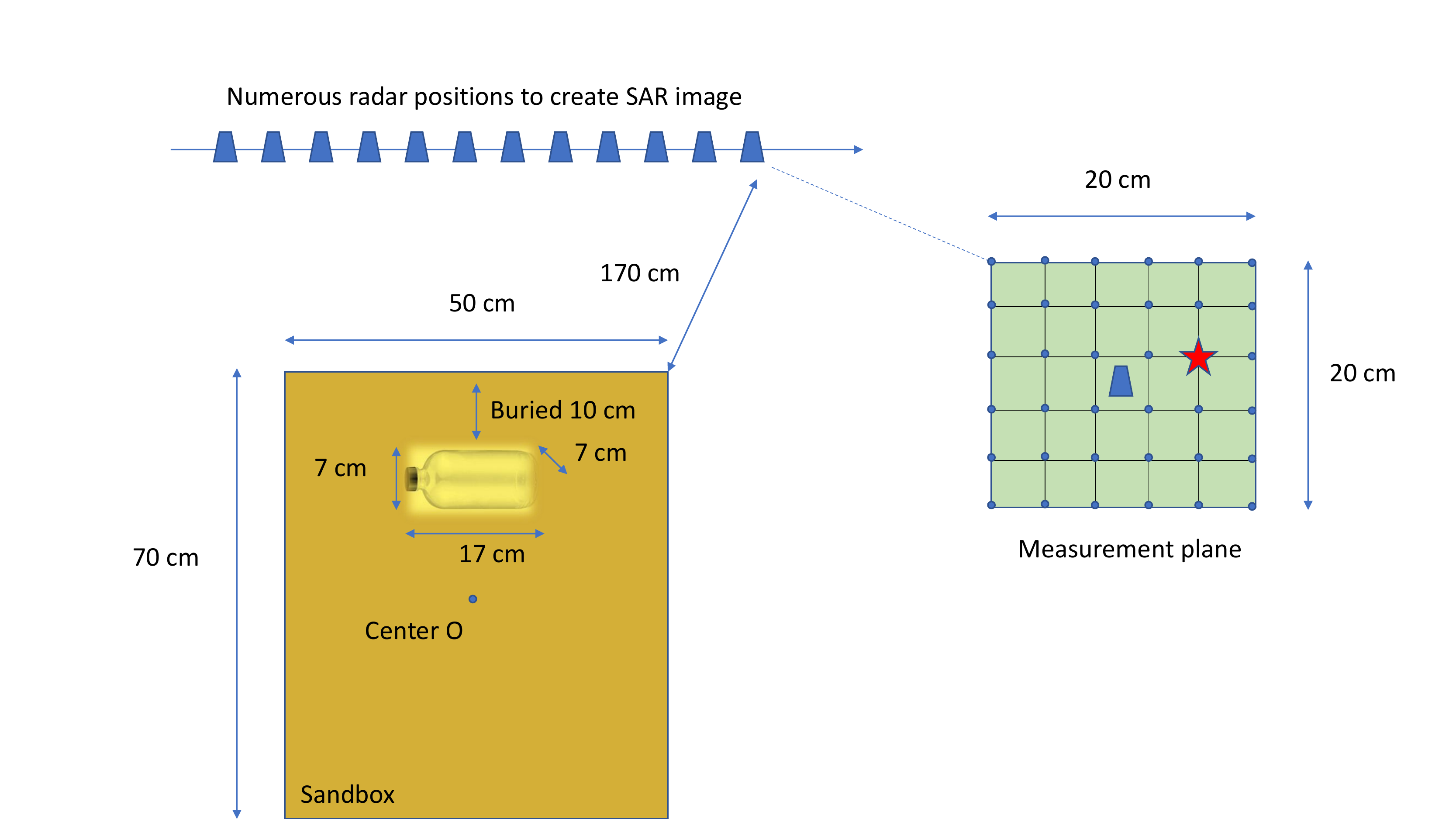}
		}
		\par\end{centering}
	\caption{(a) A schematic diagram for our
		experimental setup to collect the SAR data, which are backscattered by a buried target. (b) Basic
		geometry of the experimental data collection. In the measurement plane, the
		red star denotes the optimal detector determined in the data preprocessing
		procedure. }
	\label{fig:jug}
\end{figure}

The fake landmine under consideration was already tested in \cite{Khoa2020a,Khoa2020b} for a different mathematical model in the frequency
domain. Also different from \cite{Khoa2020a}, this time we conduct our experiment in a
room that minimizes the presence of other devices and items. This is
reasonable for both de-mining operations in the deserts and detection
services for potential threats using special drones. However, it can be seen
in Figure \ref{fig:setup} that we still have many unwanted obstacles and
furniture that can affect the quality of the backscattering data which we
measure. Mounted in front of the sandbox is a standard horn antenna with 20
cm in length, and our detector is a point probe placed behind that
transmitter. Besides, the distance between the antenna and the receiver is
about 17 cm. In order to obtain the elevation angle, the sandbox is placed
and moved diagonally to the direction that the antenna runs along with; see
Figures \ref{fig:setup} and \ref{fig:jug}. Note that due to technical
difficulties during the experiment, we actually move the sandbox instead of
running the antenna along a straight line. It still obeys the model of the
stripmap SAR imaging introduced in Figure \ref{fig:1} with a particular
elevation angle of $\pi /4$.

In our experiment, the detector is moved in a square of $20\times 20$ ($%
\text{cm}^{2}$) with the step size of 2 cm for each source position. This is
because we want the transmitters and the detectors to stay (relatively) in
the same position in the aircraft. This circumstance is different from what
we had in \cite{Khoa2020a,Khoa2020b}, where the measurement square is of $%
1\times 1$ in $\text{m}^{2}$.  Nonetheless, such a smaller working area of
the detector would hinder us from collecting good-quality signals, since the
antenna, unfortunately, is much larger than the point probe. Among many
distractions, this is the significant challenge when we work with the
experimental data in this SAR setting.

It is worth noting that our collected raw data are frequency-dependent,
while SAR imaging works with the time-dependent ones. This circumstance is
different from what we have done in \cite{Khoa2020a,Khoa2020b,Nguyen2018},
where we have tackled with the frequency domain data using the Helmholtz
equation. Therefore, we apply the inverse Fourier transform to get the
time-dependent data, which can be understood via the following operation.
Recall from section \ref{sec:41} that the original SAR data is denoted by $F(%
\mathbf{x}_{0},t)$. We then denote correspondingly the raw
frequency-dependent data by $\tilde{F}(\mathbf{x}_{0},k)$, where $k>0$ is
the wavenumber. The inverse Fourier transform reads as 
\begin{equation}\label{Fourier}
	F(\mathbf{x}_{0},t)=\dint\limits_{-\infty }^{\infty }\tilde{F}(\mathbf{x}%
	_{0},k)e^{-\text{i}kt}dk.
\end{equation}%
In practice, we actually integrate in \eqref{Fourier} over the interval $(k_1,k_2)$ defined in section \ref{sec:3.1}.

To this end, \textquotedblleft m\textquotedblright\ stands for meter;
\textquotedblleft s\textquotedblright\ and \textquotedblleft
ns\textquotedblright\ mean second and nanosecond, respectively. We remark
that the relation between the wavenumber $k$ and the frequency (denoted by $%
\text{Fre}$) is $k=2\pi \text{Fre}/c_{0}$, where $c_{0}=3\times 10^{8}$
(m/s) is the speed of light in vacuum. In this regard, we approximate the inverse
Fourier transform using the frequencies ranging from 5.6 GHz to 9.0 GHz.

As in \cite{Klibanov2021}, we compute
an approximation to the minimizer of $J_{\lambda ,\gamma }$ in \eqref{JJ}
using the standard finite difference setting. In this sense, the finite
difference operators are applied to approximate the corresponding
differential operators involved in $\mathcal{S}$; see \eqref{JJ}. For
brevity, we do not introduce those, but refer to, e.g., \cite{Khoa2020a} and
references cited therein, where we have introduced analogs of discrete
operators. Then the discrete solution is sought by the
minimization of the corresponding discrete cost functional. As in \cite%
{Klibanov2021}, we only consider the simpler penalty term in $H^{2}(R)$
instead of $H^{4}(R)$ presented in the theoretical section. This replacement
helps to reduce the complexity of computations.

It is worth mentioning that in the same vein as what we have done in \cite%
{Klibanov2021}, we scale the simulated data by a calibration factor $CF$.
This is because the problem under consideration is a quite challenging one
that requires the settings in the forward and inverse problems to be chosen
differently. We hope to address this in the future.


Consider now the speed of light $c_{0}=0.3$ (m/ns). In the case of
dimensions, the wave equation has the following form: 
\begin{equation}\label{eq:1}
	\frac{c(y)}{c_{0}^{2}}u_{ss}=u_{yy},
\end{equation}%
where we consider $s$ in ns and $y$ in m. Then let $x=y/\left( 0.3\text{m}%
\right) $ and $t=s/\text{ns}$. We find that 
\begin{equation*}
	u_{xx}=\left( 0.3\text{m}\right) ^{2}u_{yy},\quad u_{tt}=\left( 1\text{ns}%
	\right) ^{2}u_{ss}.
\end{equation*}%
Henceforth, we arrive at a new form of (\ref{eq:1}) in the dimensionless
regime: 
\begin{equation}
	\tilde{c}(x)u_{tt}=u_{xx},\quad \tilde{c}(x)=c(0.3x).  \label{eq:dimenless}
\end{equation}%
This equation resembles \eqref{forward}. In this sense, $x=1$ means $y=0.3$
m in the reality, and $t=1$ means $s=1$ ns. Thus, to construct the image
from experimental data, we solve via the convexification the CIP for
equation (\ref{eq:dimenless}) with the initial conditions as in the second line of %
\eqref{forward}, the data for the CIP as in (\ref{3.4}) and the unknown
coefficient $\tilde{c}(x).$

In our experiment, the measured data are actually far-field ones. Recall
that outside of the computational domain, which is our sandbox's length
here, the dielectric constant equals to unity. Therefore, the
dimensionless wave equation \eqref{eq:dimenless} becomes $u_{tt} = u_{xx}$
outside of that domain and up to its front surface. Suppose that the
experimental detector is positioned at $x=a$, we want to propagate the data
to a closer location $x = \tilde{a}$ with $a<\tilde{a}<0$, where the zero
point $x=0$ indicates the center of the sandbox. The last endpoint of the
sandbox is given at $x = \tilde{b}$ for $\tilde{b}>0$. According to the
absorbing boundary conditions mentioned in Remark \ref{rem:41}, the
one-dimensional wave equation satisfies the left-traveling waves, i.e. $u_t
= u_x$ at $x=a$. Similarly, the boundary condition $u_{t}=-u_{x}$ at $x=%
\tilde{b}$ is exactly satisfied by right-traveling waves.

\subsection{Data preprocessing}

In our previous works with the frequency domain data, we observed that the
signals' quality was not good; cf. e.g. \cite[Figure 3]{Nguyen2018} and \cite[Figure 3]{Khoa2020a}. That challenge is easy to understand because
of several unwanted factors coming from, e.g., the surrounding items, the
available Wi-Fi signal, etc. It was ended up with the use of the data
propagation that enables us to obtain the near-field data with a much better
quality. In the present study, we experience the same and extend the
applicability of the data propagation technique to the case of
time-dependent data. In particular, our question is to find the boundary
data $u\left( \tilde{a},t\right) $ for $a<\tilde{a}<0$ from $u\left(
a,t\right) =\varphi \left( t\right) $.

After the Fourier transform 
\begin{equation*}
	\hat{u}\left( x,k\right) =\dint\limits_{0}^{\infty }u\left(
	x,t\right) e^{\text{i}kt}dt,
\end{equation*}%
the absorbing boundary condition becomes 
\begin{equation*}
	\hat{u}_{x}\left( x,k\right) +\text{i}k\hat{u}\left( x,k\right) =0,\quad
	x\in \left[ a,\tilde{a}\right] .
\end{equation*}%
Therefore, one has 
\begin{equation*}
	\partial _{x}\left[ \hat{u}\left( x,k\right) e^{\text{i}kx}\right] =0,
\end{equation*}%
which leads to 
\begin{equation*}
	\hat{u}\left( x,k\right) e^{\text{i}kx}=\hat{u}\left( a,k\right) e^{\text{i}%
		ka}.
\end{equation*}%
Since we know the Dirichlet boundary condition measured at $x=a$, we find
that 
\begin{equation*}
	\hat{u}\left( a,k\right) =\dint\limits_{0}^{\infty }u\left(
	a,t\right) e^{\text{i}kt}dt=\hat{\varphi}\left( k\right) .
\end{equation*}%
Thus, we obtain $\hat{u}\left( x,k\right) =\hat{\varphi}\left( k\right) e^{%
	\text{i}ka}e^{-\text{i}kx}$. Then, using the inverse Fourier transform, we
derive that for any $x\in \left[ a,\tilde{a}\right] $, 
\begin{equation*}
	u\left( x,t\right) =\dint\limits_{0 }^{\infty }\hat{\varphi}\left(
	k\right) e^{-\text{i}k\left( t+x-a\right) }dk=\varphi \left( t+x-a\right) .
\end{equation*}%
Hence, $u\left( \tilde{a},t\right) =\varphi \left( t+\tilde{a}-a\right) $.
This is regarded as the time-shifted technique to move the data closer in a
spatial direction, given knowledge of the distance between the source
position and the front surface of the sandbox. Note that the operation $t+%
\tilde{a}-a$ is valid because both time and space are in the dimensionless
regime of \eqref{eq:dimenless}.

The dimensions of the sandbox and the thickness of the styrofoam give us $%
\tilde{b}=(0.35+0.07)/0.3=1.4$ and the propagated location $\tilde{a}%
=-0.5/0.3\approx -1.67$, which is a bit far away from the sandbox's front
surface covered with styrofoam. Cf. Figure \ref{fig:jug}, we can estimate
the furthest distance between the detector and the zero point centered in
the sandbox. Particularly, we compute 
\begin{align*}
	a=& -\text{maximal distance}\left( \text{antenna/sandbox }+\text{ front
		surface of sandbox/center O }\right.  \\
	& \left. +\text{ antenna/detector}\right) =-\left( 1.7+\sqrt{%
		0.42^{2}+0.25^{2}}+0.17+\sqrt{10^{2}+10^{2}}\right) /0.3 \\
	\approx & \;-8.33.
\end{align*}

Now we are in the position to complete the preprocessing procedure for our
raw far-field data.

\begin{figure}[tbp]
	\begin{centering}
		\subfloat[\label{fig:4a}]{\includegraphics[scale=0.3]{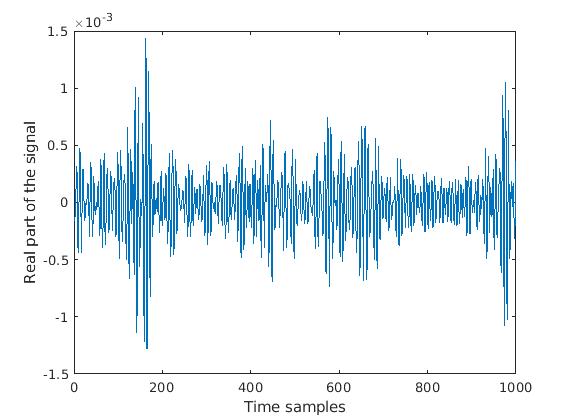}}	
		\subfloat[\label{fig:4b}]{\includegraphics[scale=0.3]{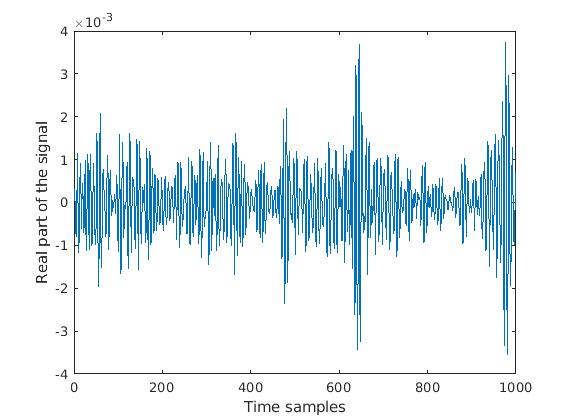}
			
		}
		\par\end{centering}
	\begin{centering}
		\subfloat[\label{fig:4c}]{\includegraphics[scale=0.3]{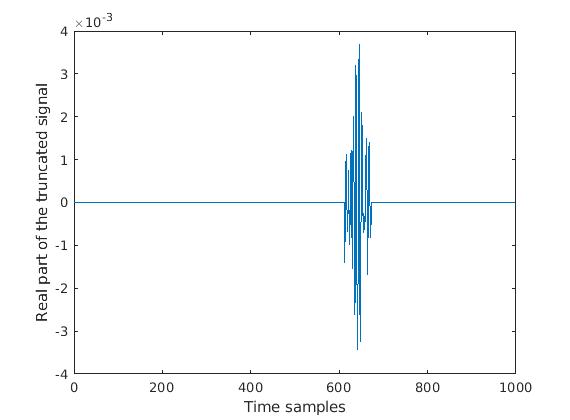}
		}
		\subfloat[\label{fig:4d}]{\includegraphics[scale=0.3]{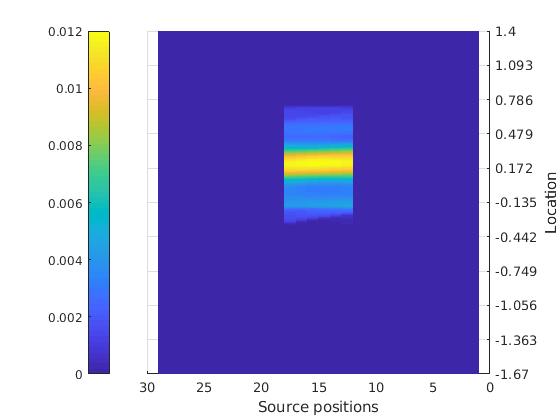}
		}
		\par\end{centering}
	\caption{A typical sample of the time-dependent experimental data for our
		experimental model after the subtraction of the background data. (a) The experimental data before data propagation. (b) The
		raw data after data propagation. (c) The corresponding experimental data
		after truncation. (d) The truncated delay-and-sum data. The delay-and-sum procedure was applied to the absolute values of the data, which is a new element here.}
	\label{fig:signal}
\end{figure}

\begin{enumerate}
	\item \textbf{Background subtraction:} For every location of the source, we subtract
	the reference data from the far-field measured data. The reference data are
	those we measure without the presence of the experimental object in the
	sandbox. It is necessary to use this subtraction since the dielectric
	constant of dry sand is around 4; cf. \cite{Daniels1996}. Currently, we
	measure the data twice since we do not know how to deal with the signal
	distractions. The subtraction is proposed to hope that it can reduce
	uncertain noises involved in our raw data. This is also another technical
	difficulty we face when working on the experiment. A similar procedure was
	implemented in \cite{Khoa2020a,Klibanov2019}. We hope to measure only once
	in the future. But to do this, we need additional analytical and numerical
	efforts.
	
	\item \textbf{Data propagation:} We apply the technique to $u(a,t)=\varphi(t)$ and
	obtain the near-field data $u\left(\tilde{a},t\right)$ for each source
	location.
	
	\item \textbf{Optimal detector:} For each source position, our detector records
	signal in a relatively small area. We found computationally that the optimal detector is determined by the maximal value of $\left\vert u\left( \tilde{a},0\right)
	\right\vert $ with respect to all detector locations. Thus, for each source
	position, we use only the data from this detector. Naturally, optimal
	detectors are different for different positions of the source. Computationally, for each source position, we note that $\left\vert u\left( \tilde{a},0\right)
	\right\vert \ne 0$ but small, which is different from the theory that it should be 0, see Figure \ref{fig:4b}. The data after this step can be denoted by $u_{\text{op}}\left(\tilde{a},t\right)$ for each source location.
	
	\item \textbf{Data truncation:} After choosing the optimal detector for each source
	position, a truncation is applied to eliminate useless signals. For every
	time-dependent signal, we can determine the location (in time) of the
	maximal absolute value of the chosen signal. In other words, such a location is found by $\max_{t\in [0,T]}\left\{\left|u_{\text{op}}\left(\tilde{a},t\right)\right|\right\}$. Then all signals outside of the interval centered around that location (with radius being 30 time samples) are
	deleted; see Figures \ref{fig:4a}, \ref{fig:4b} and \ref{fig:4c} for a
	typical sample of near-field data for our reference model. The data after this step can be denoted by $u_{\text{op}}^{\text{tr}}\left(\tilde{a},t\right)$.
	
	\item \textbf{The delay-and-sum procedure:} The so-called delay-and-sum procedure is commonly used in SAR imaging to preprocess the data \cite{Gilman2017}. We apply this procedure to the absolute value of
	the truncated data $u_{\text{op}}^{\text{tr}}\left(\tilde{a},t\right)$ and obtain $u_{\text{sad}}\left(x_{s},y_{s}\right)$ for $x_s$ being the source variable and $y_{s}$ being the location variable in the slant range plane. Then we truncate those values which are less than 95\%
	of the maximal value of the delay-and-sum data. In other words, we obtain $u_{\text{sad}}^{\text{tr}}(x_s,y_s)$ by
	\begin{align}\label{trunc}
		u_{\text{sad}}^{\text{tr}}\left(x_{s},y_s\right)=\begin{cases}
			u_{\text{sad}}\left(x_{s},y_s\right) & \text{if }u_{\text{sad}}\left(x_{s},y_s\right)\ge0.95\max_{x_{s},y_s}\left(\left|u_{\text{sad}}\right|\right),\\
			0 & \text{elsewhere}.
		\end{cases}
	\end{align}
	The preprocessed data $u_{\text{sad}}^{\text{tr}}$ are
	used then in our inversion procedure. See Figure \ref{fig:4d} for the
	processed delay-and-sum data for our reference model.
\end{enumerate}

\begin{figure}[tbp]
	\begin{centering}
		\subfloat[\label{fig:sad1}]{\includegraphics[scale=0.37]{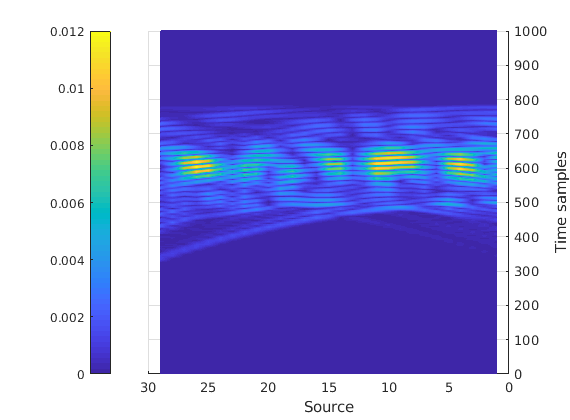}}
		\subfloat[\label{fig:sad2}]{\includegraphics[scale=0.37]{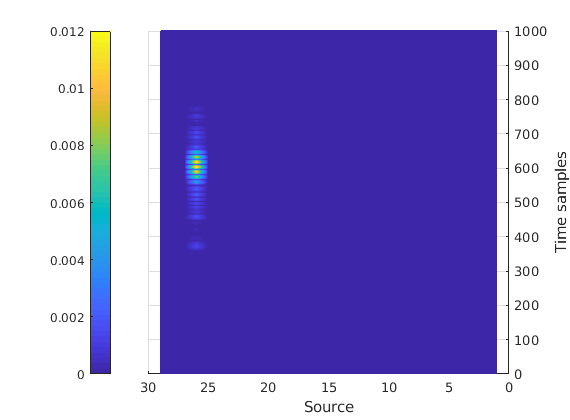}
			
		}
		\par\end{centering}
	\begin{centering}
		\subfloat[\label{fig:sad3}]{\includegraphics[scale=0.37]{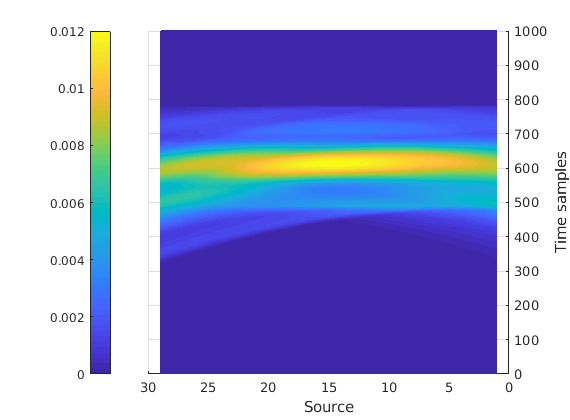}
		}
		\subfloat[\label{fig:sad4}]{\includegraphics[scale=0.37]{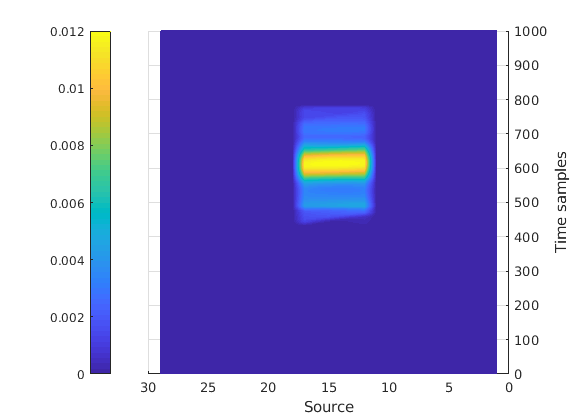}
		}
		\par\end{centering}
	\caption{(a) The delay-and-sum data when using the real part of the raw data. (b) The truncated delay-and-sum data when using the real part of the raw data. (c) The delay-and-sum data when using the absolute value of the raw data. (d) The truncated delay-and-sum data when using the absolute value of the raw data.}
	\label{fig:compare_sad}
\end{figure}

As to the data preprocessing mentioned above, we would like to remark that the conventional SAR imaging software
normally applies to two procedures to the data $F(\mathbf{x_{0}},t)$ before
getting $f_{0,n}$ for our inversion. That procedure includes the use of the
matched filtering and of the delay-and-sum technique. However, our
numerical experience in \cite{Klibanov2021} was that after the delay-and-sum alone, it provides a sufficiently good information about the
location on the slant range plane of the target represented by the slant
range image. This observation prompts one to just use the delay-and-sum technique and apply the resulting data to solve the 1D CIP along the central line of the antenna. Here, we do not detail the internal steps of the
delay-and-sum procedure and refer instead to section III-A of our first work 
\cite{Klibanov2021}. We note, however, that, compared with \cite%
{Klibanov2021}, a significantly new element of the delay-and-sum procedure
of this paper is that we apply it to absolute values of our (experimental) preprocessed data, rather
than to just real parts of the data in all previous publications. The primary reason is that the experimental data are very challenging, compared with the simulated data. By testing with the reference model, Figure \ref{fig:compare_sad} is a particular showcase of why we should use the absolute value of the data. In particular, we depict in Figures \ref{fig:sad1} and \ref{fig:sad3} the illustrations of the delay-and-sum data when using, respectively, the real value and the absolute value of the preprocessed raw data. Due to the dimensions of the sandbox and the increment of the moving antenna (see Figures \ref{fig:setup}--\ref{fig:jug} and their captions), we know that the sandbox should be between the 10-th and 20-th points of the source position. Then after the truncation \eqref{trunc}, it is clear that in terms of cleaning unnecessary data, the use of the absolute value of the data works rather better than its real value; see Figures \ref{fig:sad2} and \ref{fig:sad4}. Henceforth, if one does not have a decent data preprocessing prior solving the 1-D coefficient inverse problem, it is practically hard to obtain a good reconstruction. The rest
of formulas of the sum-and-delay procedure are the same as in section III-A of \cite{Klibanov2021}.

\subsection{Some details}

For our experimental result, we choose the calibration factor $CF=43.17$ and parameters $\lambda =2.25,\beta
=0.33,\gamma =10^{-10}$. Even
though our theoretical results are valid for large values of $\lambda $, it
is observed numerically that $\lambda \in \lbrack 1,3]$ works well. The same
interval for $\lambda $ was used in our past numerical studies; see e.g. 
\cite{Khoa2020a,Khoa2020,Klibanov2019b,Klibanov2018}. It is quite natural that the $CF$ for the
experimental data is different from the one for computationally simulated
data we did in \cite{Klibanov2021}.


\begin{figure}[tbp]
	\begin{centering}
		\subfloat[\label{fig:5a}]{\includegraphics[scale=0.065]{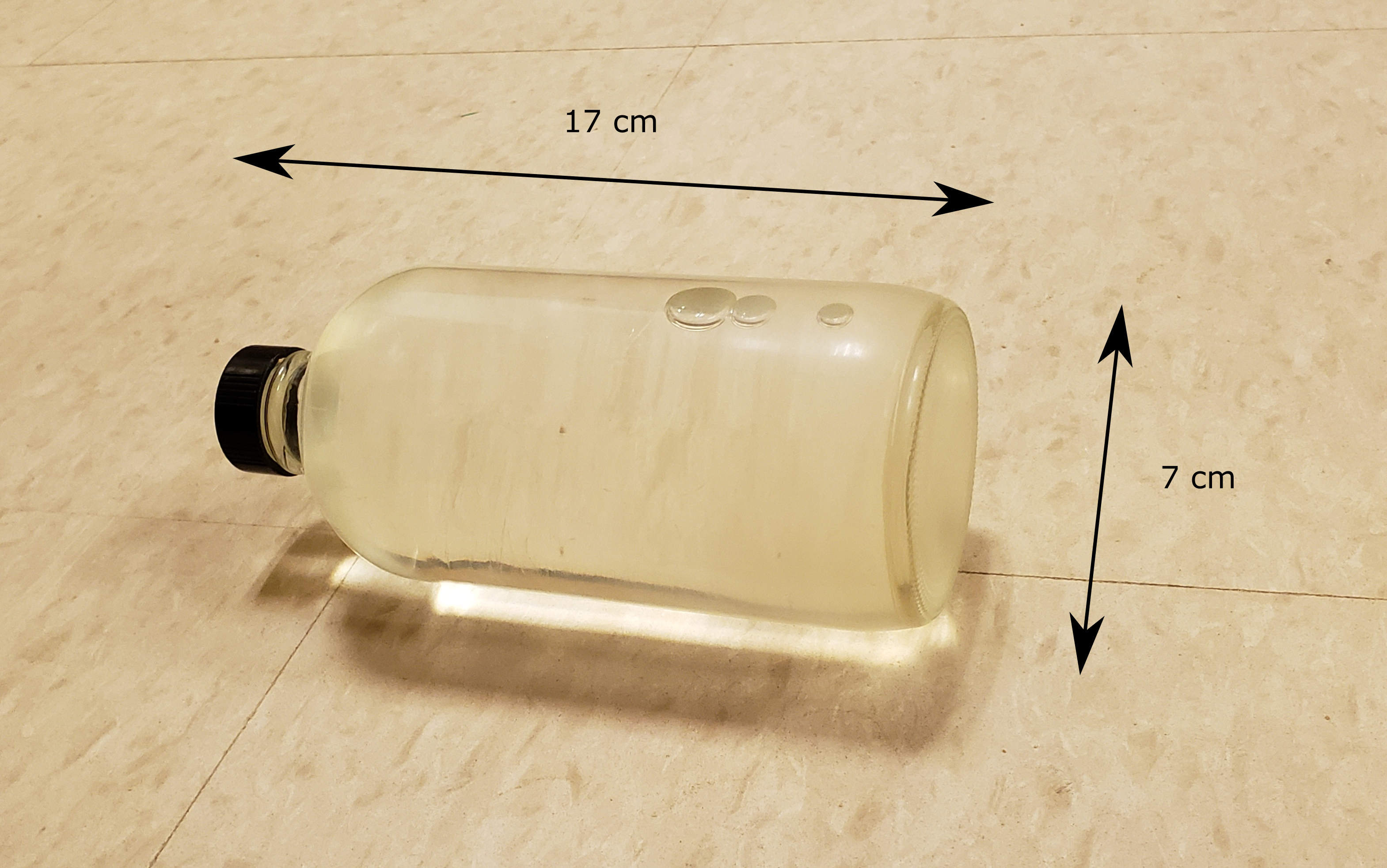}}\quad	
		\subfloat[\label{fig:5b}]{\includegraphics[scale=0.3]{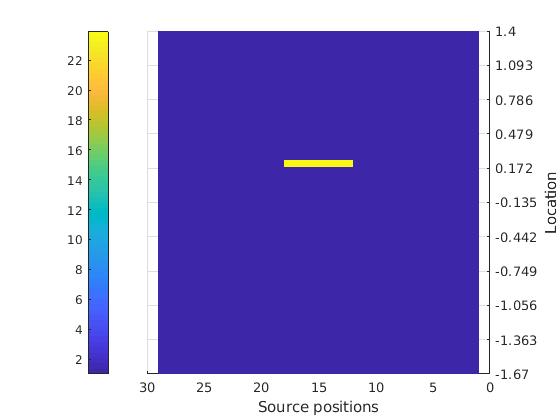}
			
		}
		\par\end{centering}
	\caption{(a) Photo of the experimental object. Its dimensions are about $17\times
		7\times 7$ in $\text{cm}^{3}$ (H$\times $W$\times $L). (b) The computed image of (a),
		where the location is presented in the dimensionless regime. The computed slant range dielectric constant is 23.9 and it was established in \protect\cite{Khoa2020a} that the true dielectric constant of this object is 23.8.}
\end{figure}

%

We remark that $T$ is determined by two times of the maximal distance
between the line of sources $L_{\text{src}}$ and the furthest point of the
sandbox. In this sense, we take $T=20>2(\tilde{b}-a)\approx 19.46$ and
consider 1000 points of time samples.

In Figure \ref{fig:signal}, we show how our data preprocessing works from
the raw data to the processed delay-and-sum data, considering the reference
model. As can be seen in Figure \ref{fig:4a}, the data propagation is
applied to avoid the direct signal from the antenna, which is as strong as
the one scattered by the target (see Figure \ref{fig:4b}). Since we do not
want any \textquotedblleft parasitic\textquotedblright\ signal, the
truncation, based on the absolute value of the signal, is used, which leads
to Figure \ref{fig:4c}.

Note again that our significant novelty here is the use of the absolute
value of the data, rather than its real part in the conventional studies of
SAR. This use is because we observe numerically that it improves the
location of the slant range target in the slant range plane, see Figure \ref{fig:compare_sad}. In the next
step, the delay-and-sum procedure with the truncation \eqref{trunc} is applied to the
resulting data, which finalizes our preprocessing procedure. Such data for
the inversion are depicted in Figure \ref{fig:4d}.

Prior to comments on the numerical results, we would like to mention a
postprocessing procedure that we use after the inversion. From the truncated
delay-and-sum data, we are able to form a rectangular area that contains
those preprocessed data of $>95$\% values of the maximal absolute value, see again \eqref{trunc}. We denote that
rectangle by $\left[ l_{1},l_{2}\right] \times \left[ s_{1},s_{2}\right] $,
where $l_{j},s_{j}$ ($j=1,2$) are standing for the location and source
coordinates, respectively. Tentatively, we denote the computed slant range
function after the inversion $\tilde{\varepsilon}_{\text{1}}(x_s,y_s)$. For each
source position $x_s$, we define $\tilde{\varepsilon}_{\text{2}}(x_s,y_s)$ as 
\begin{equation*}
	\tilde{\varepsilon}_{\text{2}}(x_s,y_s)=%
	\begin{cases}
		\max_{y_s\in \left[ l_{1},l_{2}\right] }\left( \tilde{\varepsilon}_{\text{1}%
		}(x_s,y_s)\right)  & \text{if }\tilde{\varepsilon}_{\text{1}}<0.5\max_{y_s\in \left[
			l_{1},l_{2}\right] }\left( \tilde{\varepsilon}_{\text{1}(x_s,y_s)}\right) , \\ 
		\tilde{\varepsilon}_{\text{1}(x_s,y_s)} & \text{elsewhere}.%
	\end{cases}%
\end{equation*}%
This step is to improve the dimensions of the reconstructed slant range image
since many truncation steps were applied in the data preprocessing. Then we
compute $\tilde{\varepsilon}_{\text{comp}}(x_s,y_s)$ as 
\begin{equation*}
	\tilde{\varepsilon}_{\text{comp}}(x_s,y_s)=%
	\begin{cases}
		1 & \text{outside }\left[ l_{1},l_{2}\right] \times \left[ s_{1},s_{2}\right]
		, \\ 
		\max_{\left[ l_{1},l_{2}\right] \times \left[ s_{1},s_{2}\right] }\left( 
		\tilde{\varepsilon}_{2}(x_s,y_s)\right)  & \text{inside }\left[ l_{1},l_{2}%
		\right] \times \left[ s_{1},s_{2}\right] .%
	\end{cases}%
\end{equation*}

We find that the objects' signal is approximately at the 600-th point of the
time samples. We now consider the position of the objects in terms of the distance between their central point and the central $O$ of the sandbox. Cf. Figure \ref{fig:jug3}, the central $O$ of the sandbox is 35 cm away from its front surface. As the data are propagated to $\tilde{a}=-1.67$, these time
samples range from $\tilde{a}$ to $\tilde{b}$, meaning that the distance between the front and rear surfaces of the sandbox 
should be $\tilde{b}-\tilde{a}=3.07$. Hence, the center of the object should be located $%
0.6\times (\tilde{b}-\tilde{a})-\tilde{a}=0.172$ away from the central point $O$
of the sandbox. Hence, multiplying this number by 0.3 to get back the
dimension regime, the center of the bottle is around 5 cm upwards of the central $O$. Note that ``upwards'' means we are going through the central point $O$ in the direction to the rear surface of the sandbox. The calculated position is not so accurate for the experimental bottle of water, as it should only be 21 cm downwards. We believe that the data truncation we
have imposed on the near-field data leads to this inaccuracy.


While the location in depth of our experimental object is somewhat
acceptable, we accentuate that its length is noticeably accurate. Since
the buried object is put horizontally, the length we mean here is actually
the height in reality. Counting the yellow part in Figure \ref{fig:5b} in
combination with the knowledge that the increment of source positions is 5
cm, it reveals that the length of the computed bottle is 35 cm, centered at
the 15-th source position. The true dimension should be 17 cm (cf. Figure \ref{fig:5a}), which is smaller than the computed one. However, we find it still
useful in the detection of land mines since having information about smaller
sizes is useful in the battlefield. Besides, the reconstruction is fully
accurate that the center of the bottle should be around the 15-th source
position.

Cf. \cite{Khoa2020a}, the true dielectric constant of the watered bottle is $\varepsilon _{\text{true}}=23.8$.  Compared with this, we report that its computed dielectric constant is $\tilde{\varepsilon}_{\text{comp}}=23.9$. At the same time, when $\lambda = 0$, we find that $\tilde{\varepsilon}_{\text{comp}}=15.9$, which tells the importance of the Carleman weight function presented in our reconstruction. Our numerical result shows that we are capable of reconstructing experimentally buried targets, whose dielectric constants are far away from the unity.

\section{Summary}

\label{sec:conclusions}

In our recent work \cite{Klibanov2021}, we have proposed, for the first
time, a globally convergent convexification numerical method for a fully
nonlinear model of SAR imaging. Some works on nonlinear inverse SAR problems
used locally convergent methods; cf. e.g. \cite{Karakus2020,Picard2005}. In
this work, we have presented more rigorous results. From the theoretical
standpoint, we have proved new Theorems \ref{thm:11} and \ref{thm:bigthm},
which claim the global convergence of the gradient method of the
minimization of our convexified Tikhonov-like functional. This is stronger
than our previous results for the gradient projection method (cf. e.g. \cite{Bakushinskii2017,Khoa2020a,Klibanov2019b,Klibanov2018,Klibanov2019a}) since
obviously, the second one is harder to implement. Furthermore, these
theorems explain, for the first time, our consistent observation of many
previous works on the convexification that the gradient descent method works
rather than the gradient projection method.

In this second work, we have presented a new numerical result for a quite challenging case of imaging of buried targets. This context has an obvious application in standoff identification of antipersonnel land mines and improvised explosive devices. The numerical result is for experimental data. We have used only one target for a logistical reason since the data collection process is quite time consuming. We have demonstrated that our proposed
method provides a quite accurate value of the dielectric constant of an explosive-like target
with an acceptable location and size.

\section*{Acknowledgments}

V. A. Khoa acknowledges Prof. Paul Sacks (Iowa, USA) for recent support of
his research career.

\bibliographystyle{siamplain}
\bibliography{references}
\end{document}